\documentclass[10pt]{amsart}
\usepackage{qtree,array,youngtab}
\usepackage{color, graphicx, comment}
\usepackage{tikz}
\usepackage{geometry}
\usepackage[T1]{fontenc}

\definecolor{darkgreen}{rgb}{0.,0.5,0.}

\DeclareMathOperator{\fix}{fix}
\DeclareMathOperator{\Invo}{Invo}
\def\co{d}

\DeclareMathOperator{\T}{T}

\newcommand{\la}{\lambda}

\allowdisplaybreaks

\def\setH{\mathcal{H}}
\def\setP{\mathcal{P}}

\def\O{\mathcal{O}}

\numberwithin{equation}{section}
\newtheorem{thm}{Theorem}[section]
\newtheorem{cor}[thm]{Corollary}

\newtheorem{con}[thm]{Conjecture}
\newtheorem{lem}[thm]{Lemma}

\newtheorem{prop}[thm]{Proposition}

\newtheorem{mainthm}{Theorem}

\DeclareMathOperator{\shape}{sh}

\catcode`\@=11
\def\c@rr@#1{\vbox{%
  \hrule height \ep@isseur%
   \hbox{\vrule width\ep@isseur\vbox to \t@ille{%
           \vfil\hbox  to \t@ille{\hfil#1\hfil}\vfil}%
            \vrule width\ep@isseur}%
      \hrule height \ep@isseur}}
\def\ytableau#1#2#3#4{\vbox{%
  \gdef\ep@isseur{#2}
   \gdef\t@ille{#1}
    \def\\##1{\c@rr@{$#3 ##1$}}
  \lineskiplimit=-30cm \baselineskip=\t@ille%
    \advance \baselineskip by \ep@isseur%
     \halign{%
      \hfil$##$\hfil&&\kern -\ep@isseur%
       \hfil$##$\hfil \crcr#4\crcr}}}%
\catcode`\@=12

\title[Distinct parts of Plancherel random partitions]{Distinct parts of Plancherel random partitions: exact formulas and complete asymptotic expansions}
\date{November 21, 2026}
\author[Guo-Niu Han]{Guo-Niu Han}
\author[Huan Xiong]{Huan Xiong$^*$}
\address{Universit\'e de Strasbourg, CNRS, IRMA UMR 7501, F-67000 Strasbourg, France} \email{guoniu.han@unistra.fr}

\address{Institute of Advanced Study of Mathematics,
   Harbin Institute of Technology, Heilongjiang 150001, China
   }  \email{huan.xiong.math@gmail.com}



\thanks{$^*$ Huan Xiong is the corresponding author.}

\begin{document}
\begin{abstract} 
Bufetov proved that, for a Plancherel random partition of $n$, the number of boxes with hook length $h$, divided by $\sqrt n$, converges in probability to $\frac{32h^2}{(4h^2-1)\pi^2}$; for $h=1$ this is the number of distinct parts. For the expectation we obtain an exact formula in terms of Catalan numbers, valid for every $n$ and $h$, an elementary proof of the leading term, and the complete asymptotic expansion to all orders, which is determined by three connection constants that we evaluate explicitly; written out up to $\O(n^{-9/4})$, it contains an oscillatory term of order $n^{-7/4}$. As a consequence we obtain the expansion of the expected sum of the hook lengths, extending a result of Schachinger. We also discuss the $\beta$-Plancherel measures. For the Gelfand measure ($\beta=1$), a bijection reduces the expected number of distinct parts to the mean number of fixed points of a random involution, and for every $\beta>0$ we prove that the expected number of distinct parts is at most $(\frac4\pi+o(1))\sqrt n$.
\end{abstract}
\maketitle

\section{Introduction}
Random integer partitions lie at the crossroads of many important fields in mathematics, such as combinatorics, number theory, probability theory, representation theory and mathematical physics \cite{Baik1999DJ,Baik2001Rains,Dolega2016Feray,Nekrasov2006Okounkov,Okounkov2005}. In the study of random integer partitions, a central problem that has drawn much attention is to determine their asymptotic distributions (for example, limit shapes and the number of parts) under several kinds of Plancherel-type measures. Various results on this topic have been obtained in the past decades by many mathematicians \cite{Baik2001Rains,Bernstein2007HR,Dolega2016Feray,Dolega2019Sniady,Ekhad2015Zeilberger,EvenZohar2022,Ivanov2006,Ivanov2002Olshanski,Komlos2020ST,Logan1977Shepp,Matsumoto2017,Meliot2010,Vershik1977Kerov,Zaleski2017Zeilberger}.

This paper is about the number of distinct parts of a random partition, equivalently the number of corners of its Young diagram, and more generally the number of boxes with a given hook length. For the Plancherel measure, the starting point is a result of Bufetov \cite{Bufetov2012}. In his proof of the Vershik--Kerov conjecture on the entropy of the Plancherel measure, he showed that the number of boxes with hook length $k$ of a Plancherel random partition of $n$, divided by $\sqrt n$, converges in probability to $\frac{32k^2}{(4k^2-1)\pi^2}$ \cite[Lemma~2.1]{Bufetov2012}. For $k=1$ these boxes are the corners, so that the number of distinct parts, divided by $\sqrt n$, converges to $\frac{32}{3\pi^2}$. Bufetov's proof uses the determinantal structure of the Poissonized Plancherel measure and its local convergence to the discrete sine process, and it gives no information beyond the leading order. For the expectation, the present paper goes considerably further. We obtain exact formulas, valid for every $n$; a new and elementary proof of the leading term, for every $k$; and the complete asymptotic expansion, to all orders, of the expected number of distinct parts, in which every constant is explicit; written out up to $\O(n^{-9/4})$, it contains an oscillatory term of order $n^{-7/4}$ (Theorems~\ref{mth:exact}--\ref{mth:hookh} in Section~\ref{ssec:main}). We also study the analogous questions for the $\beta$-Plancherel measures defined below.

The $\beta$-Plancherel measure~$\mathbb{P}_n^{(\beta)}$ for $\beta\geq0$ and a partition $\la$ with size $n$ 
is defined by (see \cite{Baik2001Rains})
$$
\mathbb{P}_n^{(\beta)} (\lambda)\,:= \,\frac{
(f_{\la})^\beta}{\sum_{|\mu|=n}(f_{\mu})^\beta}, 
$$
where $f_\la$ denotes the number of standard Young tableaux of shape $\lambda$. Let $\setP(n)$ be the set of partitions with size $n$, and $\la_{(n)}^{(\beta)}$ be a random partition distributed according to  $\mathbb{P}_n^{(\beta)}$. For a statistic $X$ on partitions, the expectation of $X(\la_{(n)}^{(\beta)})$ is
$$
\mathbb{E}[X(\la_{(n)}^{(\beta)})]:=\sum_{\la\in\setP(n)}X(\la)\,\mathbb{P}_n^{(\beta)}(\la)
=\frac{\sum_{|\la|=n}X(\la)\,(f_\la)^\beta}{\sum_{|\mu|=n}(f_\mu)^\beta}.
$$
When $\beta=0$, one recovers the uniform measure on $\setP(n)$; when $\beta=1$, the measure $\mathbb{P}_n^{(1)}$  is also called the Gelfand measure \cite{Meliot2010}; when $\beta=2$, one obtains the usual Plancherel measure $\mathbb{P}_n^{(2)}$, which arises naturally in the representation
theory of symmetric groups \cite{Borodin2000OO}. 
By the Frobenius formula for $f_\la$ (see Section~\ref{sec:discussion}), $\mathbb{P}_n^{(\beta)}$ is a discrete log-gas with Dyson index $\beta$, which motivates the analogy with the Gaussian Orthogonal, Unitary and Symplectic Ensembles (GOE, GUE, GSE) of random matrix theory for $\beta=1,2,4$ \cite{Mehta2004}. For $\beta=1$ and $\beta=2$ this analogy is supported by the Tracy--Widom limit theorems for the largest part recalled below \cite{Baik2001Rains,Baik1999DJ}; for $\beta=4$ it is only a motivation. Note that $\mathbb{P}_n^{(\beta)}$ is not the Jack measure with parameter $\alpha=2/\beta$ studied in \cite{Dolega2016Feray,Dolega2019Sniady}: already for $n=2$, $\mathbb{P}_2^{(\beta)}$ gives probability $\frac12$ to each of $(2)$ and $(1,1)$ for every $\beta$, whereas the Jack measure gives them $\frac1{\alpha+1}$ and $\frac{\alpha}{\alpha+1}$.

A central problem in the study of random partitions is to determine the limit behavior of $\la_{(n)}^{(\beta)}$ as $n$ tends to infinity. 

\subsection{The limit shape of random partitions.} 
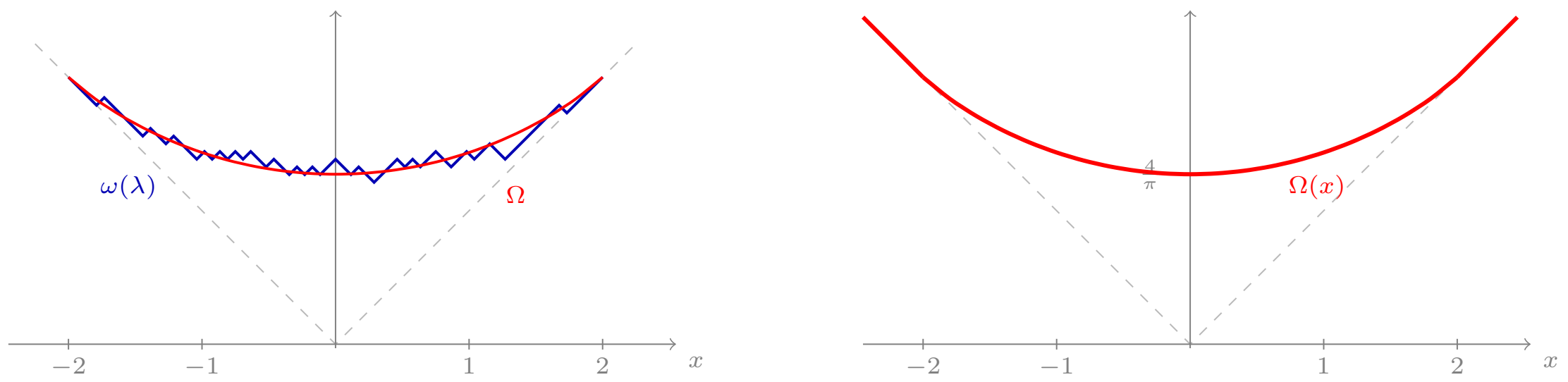
\begin{figure}[ht]
\centering
\resizebox{0.97\textwidth}{!}{%
\begin{tikzpicture}[scale=1.35]
\begin{scope}
  \draw[gray!55,dashed] (-2.25,2.25)--(0,0)--(2.25,2.25);
  \draw[->,gray] (-2.45,0)--(2.55,0) node[below right]{\scriptsize $x$};
  \draw[->,gray] (0,-0.03)--(0,2.5);
  \foreach \t in {-2,-1,1,2}
     \draw[gray] (\t,-0.04)--(\t,0.04) node[below=2pt]{\scriptsize $\t$};
  \draw[thick,blue!70!black] plot coordinates {(-2.0000,2.0000) (-1.7898,1.7898) (-1.7321,1.8475) (-1.6743,1.7898) (-1.6166,1.7321) (-1.5588,1.6743) (-1.5011,1.6166) (-1.4434,1.5588) (-1.3856,1.6166) (-1.3279,1.5588) (-1.2702,1.5011) (-1.2124,1.5588) (-1.1547,1.5011) (-1.0970,1.4434) (-1.0392,1.3856) (-0.9815,1.4434) (-0.9238,1.3856) (-0.8660,1.4434) (-0.8083,1.3856) (-0.7506,1.4434) (-0.6928,1.3856) (-0.6351,1.4434) (-0.5774,1.3856) (-0.5196,1.3279) (-0.4619,1.3856) (-0.4041,1.3279) (-0.3464,1.2702) (-0.2887,1.3279) (-0.2309,1.2702) (-0.1732,1.3279) (-0.1155,1.2702) (-0.0577,1.3279) (0.0000,1.3856) (0.0577,1.3279) (0.1155,1.2702) (0.1732,1.3279) (0.2309,1.2702) (0.2887,1.2124) (0.3464,1.2702) (0.4041,1.3279) (0.4619,1.3856) (0.5196,1.3279) (0.5774,1.3856) (0.6351,1.3279) (0.6928,1.3856) (0.7506,1.4434) (0.8083,1.3856) (0.8660,1.3279) (0.9238,1.3856) (0.9815,1.4434) (1.0392,1.3856) (1.0970,1.4434) (1.1547,1.5011) (1.2124,1.4434) (1.2702,1.3856) (1.3279,1.4434) (1.3856,1.5011) (1.4434,1.5588) (1.5011,1.6166) (1.5588,1.6743) (1.6166,1.7321) (1.6743,1.7898) (1.7321,1.7321) (2.0000,2.0000)};
  \draw[red,thick] plot[smooth] coordinates {(-2.0,2.0) (-1.8,1.8382) (-1.6,1.7085) (-1.4,1.6004) (-1.2,1.5102) (-1.0,1.4360) (-0.8,1.3765) (-0.6,1.3310) (-0.4,1.2988) (-0.2,1.2796) (0.0,1.2732) (0.2,1.2796) (0.4,1.2988) (0.6,1.3310) (0.8,1.3765) (1.0,1.4360) (1.2,1.5102) (1.4,1.6004) (1.6,1.7085) (1.8,1.8382) (2.0,2.0)};
  \node[blue!70!black] at (-1.55,1.18){\scriptsize $\omega(\lambda)$};
  \node[red] at (1.35,1.12){\scriptsize $\Omega$};
  \node at (0,-0.55){\small (a)\ \ $\omega(\lambda)$ for $|\lambda|=300$};
\end{scope}
\begin{scope}[xshift=6.4cm]
  \draw[gray!55,dashed] (-2.25,2.25)--(0,0)--(2.25,2.25);
  \draw[->,gray] (-2.45,0)--(2.55,0) node[below right]{\scriptsize $x$};
  \draw[->,gray] (0,-0.03)--(0,2.5);
  \foreach \t in {-2,-1,1,2}
     \draw[gray] (\t,-0.04)--(\t,0.04) node[below=2pt]{\scriptsize $\t$};
  \draw[densely dotted,gray] (0,1.2732)--(-0.12,1.2732) node[left]{\scriptsize $\tfrac{4}{\pi}$};
  \draw[red,very thick] plot[smooth] coordinates {(-2.0,2.0) (-1.8,1.8382) (-1.6,1.7085) (-1.4,1.6004) (-1.2,1.5102) (-1.0,1.4360) (-0.8,1.3765) (-0.6,1.3310) (-0.4,1.2988) (-0.2,1.2796) (0.0,1.2732) (0.2,1.2796) (0.4,1.2988) (0.6,1.3310) (0.8,1.3765) (1.0,1.4360) (1.2,1.5102) (1.4,1.6004) (1.6,1.7085) (1.8,1.8382) (2.0,2.0)};
  \draw[red,very thick] (-2.45,2.45)--(-2.0,2.0);
  \draw[red,very thick] (2.0,2.0)--(2.45,2.45);
  \node[red] at (0.95,1.18){\scriptsize $\Omega(x)$};
  \node at (0,-0.55){\small (b)\ \ the limit shape $\Omega$};
\end{scope}
\end{tikzpicture}%
}
\caption{(a) The rotated (Russian-convention) shape $\omega(\lambda)$ of the partition
$\lambda=(30,23^2,21,19^2,17,16,13^2,12^2,10,9,8^2,7^2,6,5,4,3^3,2^2,1^5)$ of size $300$
(exponents denote multiplicities), the Robinson--Schensted shape of a random permutation of $[300]$, together with the limit shape $\Omega$.
(b) The function $\Omega(x)$. The dashed wedge is $|x|$, to which both profiles are asymptotic;
$\Omega(0)=4/\pi$.}\label{fig:1}
\end{figure}

The limit shape of $\la_{(n)}^{(\beta)}$ has been studied since the 1970s \cite{Logan1977Shepp,Vershik1977Kerov}. 
 One can associate each partition $\la$ with a continuous piecewise linear function
$\omega(\la)$ from $\mathbb{R}$ to $\mathbb{R}$ (see \cite{Dolega2016Feray,Ivanov2002Olshanski}), which is called the shape of $\la$. More precisely, the shape $\omega(\la)$ is obtained by rotating the upper right boundary of the generalized Young diagram of $\la$ by $45^{\circ}$, and then scaling by a factor~$\sqrt{2}/\sqrt{|\la|}$. An example of the shape of a partition $\la$ is given on the left-hand side of Figure~\ref{fig:1}.

When $\beta=2$, 
the limit shape $\Omega$ of the random partition $\la_{(n)}^{(2)}$ was
discovered by Vershik--Kerov \cite{Vershik1977Kerov} and Logan--Shepp \cite{Logan1977Shepp} independently.
More precisely, they showed that 
$$
\lim_{n\to\infty}\sup_{x\in\mathbb{R}}|\, \omega\left(\la_{(n)}^{(2)}\right)(x)-\Omega(x)|=0 \quad
\text{\rm in probability},
$$ 
where (see the right-hand side of Figure~\ref{fig:1})
$$
\Omega(x)= \begin{cases} \frac{2}{\pi}(x\arcsin\tfrac x2+\sqrt{4-x^2}), \qquad -2\leq x\leq 2; \\
|x|, \qquad |x|>2.
\end{cases}
$$ 
Results on limit shapes under related Plancherel-type measures, which in general differ from $\mathbb{P}_n^{(\beta)}$, can be found in \cite{Bernstein2007HR,Dolega2016Feray,Ivanov2006,Meliot2010}.

\subsection{Number of parts.}
The distribution of the number of parts of $\la_{(n)}^{(\beta)}$ has also been widely studied. It is the same as the distribution of its largest part: since $f_{\la'}=f_\la$, where $\la'$ is the conjugate of $\la$, the measure $\mathbb{P}_n^{(\beta)}$ is invariant under conjugation.

When $\beta =0$, Kessler and Livingston \cite{Kessler1976Livingston} showed that 
$$
\mathbb{E}[\ell(\la_{(n)}^{(0)})] = \frac{\sqrt{6}}{2\pi} \sqrt{n} \log n + \O(\sqrt{n}),
$$
where $\ell(\la)$ denotes the length (i.e., the number of parts) of a partition $\la$.

When $\beta=1$, by the Robinson--Schensted--Knuth correspondence \cite{Knuth1998Vol3}, the distribution of the number of parts of $\la_{(n)}^{(1)}$ is the same as the distribution of the length of the longest increasing subsequence of a uniform random involution of $[n]$.
In 2001, Baik and Rains \cite{Baik2001Rains} derived the following Tracy--Widom limit theorem for $\ell(\la_{(n)}^{(1)})$:
\begin{equation}
  \lim_{n\to\infty} \Pr\biggl(\frac{\ell(\la_{(n)}^{(1)})-2\sqrt{n}}{n^{1/6}} \leq x\biggr) = {\T}_1(x),
\end{equation}
where $\T_1(x)$ is the limiting distribution function \cite{Tracy1996Widom} for the (scaled) largest eigenvalue of random matrices taken from the Gaussian Orthogonal Ensemble.

When $\beta=2$, it is known that
$$
\mathbb{E}[\ell(\la_{(n)}^{(2)})]\sim2\sqrt{n},
$$
where, for two sequences $\{a_n\}$ and $\{b_n\}$, $a_n\sim b_n$ means that $\lim_{n\to\infty}a_n/b_n=1$. The limit shape theorem alone does not imply this, since it does not control a thin first column. Vershik and Kerov \cite{Vershik1977Kerov} proved that $\ell(\la_{(n)}^{(2)})/\sqrt n\to2$ in probability, and the convergence of the expectation follows, for instance, from the convergence of moments proved in \cite{Baik1999DJ}.
A more precise limit theorem, with Tracy--Widom fluctuations, for $\ell(\la_{(n)}^{(2)})$ was derived in \cite{Baik1999DJ}:
\begin{equation}
  \lim_{n\to\infty} \Pr\biggl(\frac{\ell(\la_{(n)}^{(2)})-2\sqrt{n}}{n^{1/6}} \leq x\biggr)=
  {\T}_2(x),
\end{equation}
where ${\T}_2(x)$ is the Tracy--Widom distribution function \cite{Tracy1994Widom} for the (scaled) largest eigenvalue of the Gaussian Unitary Ensemble.

\subsection{The number of distinct parts of partitions under the uniform measure.}  Another concept that has attracted much attention is the number of distinct parts of a random partition.  Let $\co(\la)$ be the number of distinct parts of $\la$. When $\beta=0$, the central limit theorem for $\co(\la_{(n)}^{(0)})$ was derived by Goh and Schmutz \cite{Goh1995Schmutz}:
\begin{equation}
  \lim_{n\to\infty} \Pr\biggl(\co(\la_{(n)}^{(0)}) \leq (\sqrt{6}/\pi) \sqrt{n}+  \sqrt{  \sqrt{6}/2\pi - \sqrt{54}/\pi^3  } \cdot n^{1/4} x\biggr) = \Phi(x),
\end{equation}
where $\Phi(x)$ is the cumulative distribution function of the standard normal distribution. The properties of $\co(\la_{(n)}^{(0)})$ have also been studied by \cite{Blecher2016etal,Hirschhorn2014,Hwang1997Yeh,Wilf1983}. In particular, the generating function and several asymptotic formulas for $\sum_{\mid\la\mid=n} \co(\la)$ are given in \cite{Blecher2016etal}.  


\subsection{Main results.}\label{ssec:main}
Much less is known for $\beta>0$. For $\beta=2$, Bufetov's law of large numbers \cite[Lemma~2.1]{Bufetov2012}, recalled at the beginning of the introduction, gives $\co(\la_{(n)}^{(2)})/\sqrt n\to\frac{32}{3\pi^2}$ in probability, since the boxes with hook length $1$ are the corners (Lemma~\ref{lem:corner_hook_one}).

In this paper we study the expectation of $\co(\la_{(n)}^{(\beta)})$ and, for $\beta=2$, of the number $\co_h(\la)$ of boxes of $\la$ with hook length $h$, so that $\co_1=\co$. The emphasis is on exact formulas and, for $\beta=2$, on the complete asymptotic expansion. Write
$$
F_n^{(h)}:=\mathbb{E}[\co_h(\la_{(n)}^{(2)})],\qquad F_n:=F_n^{(1)}=\mathbb{E}[\co(\la_{(n)}^{(2)})],
$$
and let $C_k=\frac{(2k)!}{k!\,(k+1)!}$ be the Catalan numbers.

\begin{mainthm}[Exact formulas]\label{mth:exact}
For all integers $h\geq1$ and $n\geq1$,
$$
F_n^{(h)}=\sum_{m=h}^{n}(-1)^{m-h}\binom{n}{m}\,
\frac{h^2\,(m+1)!}{m\,(m-h)!\,(m+h)!}\,C_{m-1}C_m;
$$
in particular,
$$
F_n=\sum_{m=1}^{n}(-1)^{m-1}\binom nm\frac{C_{m-1}C_m}{m!}.
$$
Moreover, $F_n^{(h)}$ satisfies a homogeneous linear recurrence of order $5$ whose coefficients are polynomials in $n$ and $h^2$, and $F_n$ satisfies an inhomogeneous linear recurrence of order $4$.
\end{mainthm}

The exact formula is proved in Theorem~\ref{thm:exacth}, and the recurrences in Propositions~\ref{prop:rec} and~\ref{prop:rechall}. The alternating sums of Theorem~\ref{mth:exact} hide the growth of $F_n^{(h)}$: their terms are much larger than their values. Our main result is the complete asymptotic expansion of $F_n$, with all constants explicit.

\begin{mainthm}[Distinct parts]\label{mth:full}
\textup{(a)} There are real constants $K_1,K_2$ and a complex constant $K_3$ such that, for $n$ large enough,
\begin{equation}\label{eq:fullexp}
F_n=K_1\,S_1(n)+K_2\,S_2(n)-\frac12+K_3\,Y(n)+\bar K_3\,\bar Y(n),
\end{equation}
where $S_1,S_2,Y$ are the solutions of the homogeneous recurrence of order $4$ satisfied by $F_n+\frac12$ whose asymptotic expansions are the formal solutions of Proposition~\ref{prop:formal}, and $\bar Y$ is the complex conjugate of $Y$.

\textup{(b)} The constants are
$$
K_1=\frac{32}{3\pi^2},\qquad
K_2=\frac{4}{3\pi^2}+\frac{2\gamma+12\log2-1}{8\pi^2},\qquad
K_3=-\frac{e^{8+i\pi/4}}{1024\,\pi^{3/2}},
$$
where $\gamma$ is Euler's constant.
\end{mainthm}

Part~(a) is proved in Section~\ref{ssec:proofshape}. Part~(b) is proved in Section~\ref{ssec:distinct} for $K_1$, in Section~\ref{ssec:k2evaluation} for $K_2$, and in Section~\ref{ssec:cconst} for $K_3$. Together with the coefficients of Proposition~\ref{prop:formal}, Theorem~\ref{mth:full} gives
\begin{equation}\label{eq:Fexplicit}
\begin{aligned}
F_n={}&\frac{32}{3\pi^2}\sqrt n-\frac12+\frac{3\log n+6\gamma+36\log2+29}{24\pi^2\sqrt n}
-\frac{1146\log n+2292\gamma+13752\log2-6251}{3\cdot2^{13}\pi^2\,n^{3/2}}\\
&-\frac{e^{8}\cos\bigl(8\sqrt n+\frac\pi4\bigr)}{512\,\pi^{3/2}\,n^{7/4}}+\O\bigl(n^{-9/4}\bigr).
\end{aligned}
\end{equation}
Bufetov's theorem gives the leading term $K_1\sqrt n$ (for the expectation, after a uniform integrability argument), and nothing beyond it. Our proof of the leading term is different and elementary, and all the lower-order terms, including the constants $K_2$ and $K_3$, are new. The same method applies to every hook length.

\begin{mainthm}[Every hook length]\label{mth:hookh}
Let $h\geq1$ be a fixed integer; in the asymptotic statements below, $n\to\infty$ with $h$ fixed.

\textup{(a)} There are real constants $K_1^{(h)},K_2^{(h)}$ and a complex constant $K_3^{(h)}$ such that, for $n$ large enough \textup(depending on $h$\textup),
\begin{equation}\label{eq:fullexph}
F_n^{(h)}=K_1^{(h)}S_1^{(h)}(n)+K_2^{(h)}S_2^{(h)}(n)-\frac h2+K_3^{(h)}Y^{(h)}(n)+\bar K_3^{(h)}\,\overline{Y^{(h)}}(n),
\end{equation}
where $S_1^{(h)},S_2^{(h)},Y^{(h)}$ are the solutions of the recurrence of order $5$ of Theorem~\ref{mth:exact} whose asymptotic expansions are the formal solutions of Proposition~\ref{prop:formalh}, and $\overline{Y^{(h)}}$ is the complex conjugate of $Y^{(h)}$.

\textup{(b)} The constants are
$$
K_1^{(h)}=\frac{32h^2}{(4h^2-1)\pi^2},\qquad
K_3^{(h)}=(-1)^{h-1}h^2K_3,
$$
$$
K_2^{(h)}=\frac{h^2}{\pi^2}\Bigl(\frac{4}{4h^2-1}+\frac{2\gamma+12\log2+3}{8}-\frac12\sum_{i=1}^{h}\frac1{2i-1}\Bigr).
$$
\end{mainthm}

For $h=1$, Theorem~\ref{mth:hookh} is Theorem~\ref{mth:full}. Part~(a) is proved in Section~\ref{ssec:fullh}. Part~(b) is proved in Section~\ref{ssec:asyh} for $K_1^{(h)}$, and in Section~\ref{ssec:consth} for $K_2^{(h)}$ and $K_3^{(h)}$. Together with the coefficients of Proposition~\ref{prop:formalh}, Theorem~\ref{mth:hookh} gives
\begin{equation}\label{eq:Fhexplicit}
\begin{aligned}
F_n^{(h)}={}&\frac{32h^2}{(4h^2-1)\pi^2}\sqrt n-\frac h2
+\frac{h^2}{8\pi^2\sqrt n}\Bigl(\log n+2\gamma+12\log2+3+\frac{32}{4h^2-1}-4\sum_{i=1}^{h}\frac1{2i-1}\Bigr)\\
&+\frac{h^2(4h^2-577)}{3\cdot2^{12}\pi^2\,n^{3/2}}\Bigl(\log n+2\gamma+12\log2+3-4\sum_{i=1}^{h}\frac1{2i-1}\Bigr)\\
&+\frac{h^2(-368h^4+37144h^2+9169)}{9\cdot2^{13}(4h^2-1)\pi^2\,n^{3/2}}+(-1)^h\frac{h^2e^{8}\cos\bigl(8\sqrt n+\frac\pi4\bigr)}{512\,\pi^{3/2}\,n^{7/4}}+\O_h\bigl(n^{-9/4}\bigr),
\end{aligned}
\end{equation}
where the constant implied in $\O_h$ depends on $h$. By contrast, the exact formula of Theorem~\ref{mth:exact} holds for all $n$ and $h$.
For the expectation, the leading term of Theorem~\ref{mth:hookh} is Bufetov's constant $\frac{32h^2}{(4h^2-1)\pi^2}$; our proof uses only the exact formula of Theorem~\ref{mth:exact}.

Schachinger \cite{Schachinger2023} obtained, by Rice's integral, the asymptotic expansion of the expected sum of the hook lengths $\mathbb{E}\bigl[\sum_{\square\in\la}h_\square\bigr]$ up to $\O(n^{-9/4})$. It involves the same combination $\log n+2\gamma+12\log2$ and the same oscillation $\cos(8\sqrt n+\frac\pi4)\,n^{-7/4}$ as \eqref{eq:Fexplicit}. The reason is an exact identity between the two expectations (Proposition~\ref{prop:hsidentity}). As a consequence of Theorem~\ref{mth:full}, we obtain the complete asymptotic expansion of $\mathbb{E}\bigl[\sum_{\square}h_\square\bigr]$, with all constants explicit (Corollary~\ref{cor:hooksum}); this gives in particular an independent proof of his result.

For $\beta=0$ and $\beta=1$, the complete asymptotic expansion of $\mathbb{E}[\co(\la_{(n)}^{(\beta)})]$ is already known. For the uniform measure ($\beta=0$), it was obtained by Grabner, Knopfmacher and Wagner \cite[Proposition~4.1]{Grabner2014KW}:
$\mathbb{E}[\co(\la_{(n)}^{(0)})]=\frac{\sqrt6}{\pi}\sqrt n+\frac3{\pi^2}-\frac12+\O(n^{-1/2})$. For the Gelfand measure ($\beta=1$), a bijection based on reverse row insertion (Theorem~\ref{th:fix=cor_sec1}) shows that $\sum_{|\la|=n}f_\la\,\co(\la)$ is the total number of fixed points of the involutions of $[n]$. Hence $\mathbb{E}[\co(\la_{(n)}^{(1)})]$ is the mean number of fixed points of a uniform random involution, whose complete expansion was given by Lugo \cite[p.~188]{Lugo2010}: $\mathbb{E}[\co(\la_{(n)}^{(1)})]=\sqrt n-\frac12+\frac{3}{8\sqrt n}+\O(n^{-1})$. See Section~\ref{ssec:beta01}. For general $\beta>0$ we prove an upper bound.

\begin{mainthm}[General $\beta$]\label{mth:beta}
For every $\beta>0$,
$\displaystyle\limsup_{n\to\infty}\mathbb{E}[\co(\la_{(n)}^{(\beta)})]/\sqrt n\le\frac4\pi$.
\end{mainthm}

Theorem~\ref{mth:beta} is proved in Section~\ref{ssec:cupper}. The leading constants $c(0)=\frac{\sqrt6}{\pi}$, $c(1)=1$ and $c(2)=\frac{32}{3\pi^2}$ increase with $\beta$, and we conjecture that $\mathbb{E}[\co(\la_{(n)}^{(\beta)})]\sim c(\beta)\sqrt n$ with $c(\beta)$ non-decreasing (Conjecture~\ref{con:cbeta}). In Section~\ref{sec:discussion} we also discuss the limits of $c(\beta)$ as $\beta\to\infty$ and $\beta\to0^+$ (Conjecture~\ref{con:cinfty}).

\subsection{Organization.}
Section~\ref{sec:prelim} recalls basic definitions. Sections~\ref{sec:exact}--\ref{sec:constants} treat the Plancherel measure: the exact formulas, the recurrences and the integral representations (Section~\ref{sec:exact}), the leading asymptotics (Section~\ref{sec:beta2}), the shape of the complete expansion (Section~\ref{sec:rec}) and the evaluation of the connection constants (Section~\ref{sec:constants}), which ends with the expected sum of the hook lengths (Section~\ref{ssec:hooksum}). Section~\ref{sec:discussion} discusses the $\beta$-Plancherel measures: the known results for $\beta=0$ and $\beta=1$, together with the fixed-point bijection for $\beta=1$, the upper bound of Theorem~\ref{mth:beta}, and two conjectures on~$c(\beta)$.

\section{Preliminaries}\label{sec:prelim}
In this section, we fix some notation and give some preliminary results. 

The basic knowledge on partitions and permutations can be found in \cite{Knuth1998Vol3,Macdonald1995,Stanley1999EC2}.
A \emph{partition} is a finite non-increasing sequence of  positive
integers $\lambda = (\la_1, \la_2, \ldots, \la_\ell)$, where $\ell(\la):=\ell$ is called the length, 
 $\la_i \ (1\leq i\leq \ell)$ are called the parts, and $|\lambda|=\sum_{1\leq i\leq \ell}\la_i$ is the size of $\la$.
Recall that $\co(\la)$ is the number of distinct parts of $\la$. Let $\setP(n)$ be the set of partitions with size $n$ and $\setP$ be the set of all partitions (including the empty partition).
A partition~$\lambda$ can be visually represented by its Young diagram $Y(\la)$, which is a collection of boxes arranged in a left-justified way such that there are exact $\la_i$ boxes in the $i$-th row for $1\leq i\leq \ell(\la)$. 

A box is called a corner if there are no other boxes exactly below or to the right of it. For each box $\square$ in the Young diagram $Y(\la)$, the \emph{hook length} $h_{\square}$ of $\square$ is the number of boxes exactly below, exactly to the right, and the box itself. For example,
the Young diagram of the partition $(4,3,3,2)$ and its hook lengths are given in Figure \ref{fig:2.1}. Let $\setH(\la)$ be the multiset of all hook lengths of $\lambda$ (each box contributing one element) and $H_{\lambda}=\prod_{\square\in\lambda} h_{\square}$ be the product
of all hook lengths of $\lambda$.  A standard Young tableau (SYT) of shape $\la$ is a filling in all
boxes of $Y(\la)$ with distinct numbers from $1$ to $|\la|$ such
that the numbers are increasing in each row and each column.
The number of standard Young tableaux of shape
$\lambda$ is denoted by $f_\lambda$. The following relation between $f_\lambda$ and  $H_{\lambda}$ was first discovered by Frame, Robinson and Thrall \cite{Frame1954RT}:
\begin{align}
f_\lambda = \frac{|\lambda| !}{H_{\lambda}}.
\end{align}

\begin{figure}[htbp]
\begin{center}
\Yvcentermath1

\begin{tabular}{c}

$\yng(4,3,3,2)$
\ \ \ \ \ \ \ \
$\young(7641,542,431,21)$

\end{tabular}

\end{center}
\caption{The Young diagram of the partition $(4,3,3,2)$ and its hook
lengths.}
\label{fig:2.1}
\end{figure}

A permutation $\sigma$ of length $n$ is defined as a bijection from $\{1,2,3,\cdots,n\}$ to itself. A permutation $\sigma$ is called an involution if $\sigma^2$ is the identity map of the set $\{1,2,3,\cdots,n\}$.

The following lemma shows that the number of distinct parts for a partition of size $n$ is at most~$\sqrt{2n}$.  
\begin{lem}\label{lem: upper_bound}
Let $\la$ be a partition with size $n\geq 1$. Then 
$$
\co(\la) < \sqrt{2n}.
$$
\end{lem}
\begin{proof}
Suppose that $\la$ has $k=\co(\la)$ distinct parts. Then these $k$ distinct parts are at least $1,2,\ldots,k$. Thus
$$
n=|\la|\geq 1+2+\cdots+k=\frac{k(k+1)}{2}>\frac{k^2}{2}.
$$
Therefore, 
$$
k<\sqrt{2n}.
$$
\end{proof}

It is easy to see that the number of distinct parts in a partition $\la$ is equal to the number of corners in the corresponding Young diagram, thus is equal to the number of boxes with hook length $1$ in the Young diagram. Therefore, the following result holds.

\begin{lem}\label{lem:corner_hook_one}
Let $\la$ be a partition. Then we have
$$
\co(\la)= \#\{\square\in\la, h_{\square}=1\}.
$$
\end{lem}

For example, by Figure \ref{fig:2.1} we can see the number of boxes with hook length $1$ in the Young diagram of the partition $(4,3,3,2)$ is equal to $3$, which is equal to the number of distinct parts of the partition $(4,3,3,2)$.

\section{Exact formulas under the Plancherel measure}\label{sec:exact}

In Sections~\ref{sec:exact}--\ref{sec:constants} we consider the Plancherel measure. Write
$$
F_n:=\mathbb{E}[\co(\la_{(n)}^{(2)})]=\frac{1}{n!}\sum_{|\la|=n}(f_\la)^2\,\co(\la),
$$
and let $\mathbb{E}$ denote the expectation with respect to the Plancherel measure $\mathbb{P}_n^{(2)}$, so that $\mathbb{P}_n^{(2)}(\la)=(f_\la)^2/n!$.

\subsection{Distinct parts.}\label{ssec:exact1}
We first prove the exact formula for $F_n$.

\begin{thm}\label{th:exactF2_sec1}
We have
	\begin{equation}
		\mathbb{E}[\co(\la_{(n)}^{(2)})]
		=\sum_{k=0}^{n-1} \frac { (-1)^k C_k C_{k+1} }{(k+1)!} \binom {n}{k+1},
	\end{equation}
	where
	$$
	C_k= \frac {(2k)!}{k!(k+1)!}
	$$
	is the $k$-th Catalan number.
\end{thm}

For integers $r\geq 1$ consider the Plancherel averages of hook-length products
$$
A_r(n):=\mathbb{E}\Bigl[\sum_{\square\in\la}\prod_{1\le j\le r}(h_\square^2-j^2)\Bigr],\qquad
B_r(n):=\mathbb{E}\Bigl[\sum_{\square\in\la}\prod_{2\le j\le r}(h_\square^2-j^2)\Bigr].
$$
The starting point is the following remarkable evaluation. Nekrasov and Okounkov \cite{Nekrasov2006Okounkov} discovered the hook-length formula
$$
\sum_{n\ge0}\frac{x^n}{(n!)^2}\sum_{|\la|=n}(f_\la)^2\prod_{\square\in\la}(t+h_\square^2)
=\prod_{i\ge1}(1-x^i)^{-1-t};
$$
an elementary proof and a refinement were given by the first author \cite{Han2010Fourier}, who conjectured \cite{Han2009Conj} that the Plancherel average of $\sum_{\square}h_\square^{2k}$ is a polynomial in $n$ for every $k$. This was proved by Stanley \cite{Stanley2010}, which led Okada to conjecture the explicit polynomial below; it was proved by Panova \cite{Panova2012} (see also \cite{Han2018Xiong} for difference-operator proofs and extensions).

\begin{lem}[Okada--Panova \cite{Panova2012}]\label{lem:okada}
For all integers $r\ge1$ and $n\ge0$,
$$
A_r(n)=\Lambda_r\binom{n}{r+1},\qquad \Lambda_r=\frac{(2r)!\,(2r+1)!}{r!\,(r+1)!^2}.
$$
\end{lem}

\begin{lem}\label{lem:Bnn}
For $n\geq 2$ we have
$\displaystyle B_n(n)=\Bigl(\prod_{j=2}^n(1-j^2)\Bigr)\,F_n.$
\end{lem}

\begin{proof}
Since $1\le h_\square\le n$, the product $\prod_{2\le j\le n}(h_\square^2-j^2)$ vanishes
unless $h_\square=1$, in which case it equals $\prod_{j=2}^n(1-j^2)$. By Lemma~\ref{lem:corner_hook_one}
the number of boxes with $h_\square=1$ is $\co(\la)$, so
$\sum_\square\prod_{2\le j\le n}(h_\square^2-j^2)=\prod_{j=2}^n(1-j^2)\,\co(\la)$;
taking the expectation gives the claim.
\end{proof}

\begin{lem}\label{lem:Brec}
We have $B_1(n)=n$, $A_1(n)=\tfrac{3n^2-3n}2$, and for $r\ge2$
\begin{equation}\label{eq:Brec}
B_r(n)=A_{r-1}(n)+(1-r^2)B_{r-1}(n);
\end{equation}
in particular $B_2(n)=\tfrac{3n^2-9n}2$. Consequently, for $r\ge2$,
\begin{equation}\label{eq:Bexpand}
B_r(n)=\sum_{k=2}^{r-1}A_k(n)\prod_{j=k+2}^{r}(1-j^2)+B_2(n)\prod_{j=3}^{r}(1-j^2).
\end{equation}
\end{lem}

\begin{proof}
We have $B_1(n)=\mathbb{E}\bigl[\sum_{\square\in\la}1\bigr]=n$, and $A_1(n)=\Lambda_1\binom n2=3\binom n2$ by Lemma~\ref{lem:okada}. For $r\ge2$ and every box, with $h=h_\square$,
$$
\prod_{2\le j\le r}(h^2-j^2)-\prod_{1\le j\le r-1}(h^2-j^2)
=\prod_{2\le j\le r-1}(h^2-j^2)\bigl[(h^2-r^2)-(h^2-1)\bigr]
=(1-r^2)\!\!\prod_{2\le j\le r-1}\!\!(h^2-j^2).
$$
Summing over the boxes and taking the expectation gives $B_r(n)-A_{r-1}(n)=(1-r^2)B_{r-1}(n)$, whence $B_2=A_1-3B_1=\tfrac{3n^2-9n}2$. Iterating
\eqref{eq:Brec} yields \eqref{eq:Bexpand}.
\end{proof}

\begin{proof}[Proof of Theorem~\ref{th:exactF2_sec1}]
For $n=1$ both sides are equal to $1$; let $n\geq 2$.
Take $r=n$ in \eqref{eq:Bexpand} and divide by $\prod_{j=2}^n(1-j^2)$. Since
$\prod_{j=3}^n(1-j^2)/\prod_{j=2}^n(1-j^2)=1/(1-2^2)=-\tfrac13$, Lemma~\ref{lem:Bnn} gives
$$
F_n=\frac{B_n(n)}{\prod_{j=2}^n(1-j^2)}
=\sum_{k=2}^{n-1}A_k(n)\prod_{j=2}^{k+1}(1-j^2)^{-1}-\tfrac13 B_2(n).
$$
Now $\prod_{j=2}^{k+1}(1-j^2)=\prod_{j=2}^{k+1}\bigl(-(j-1)(j+1)\bigr)=(-1)^k k!\,(k+2)!/2$, and $A_k(n)=\Lambda_k\binom{n}{k+1}$ by Lemma~\ref{lem:okada}. Since
$2\,(2k)!\,(2k+1)!\,(k+1)!^3=(2k)!\,(2k+2)!\,k!\,(k+1)!^2$, we obtain, for
$k\ge2$,
$$
A_k(n)\prod_{j=2}^{k+1}(1-j^2)^{-1}
=\frac{2(-1)^k(2k)!\,(2k+1)!}{k!^2\,(k+1)!^2\,(k+2)!}\binom{n}{k+1}
=\frac{(-1)^k C_kC_{k+1}}{(k+1)!}\binom{n}{k+1}.
$$
Finally $-\tfrac13 B_2(n)=\tfrac{3n-n^2}2=n-\binom n2$ is exactly the sum of the terms $k=0$ and $k=1$
of the right-hand side of Theorem~\ref{th:exactF2_sec1}, so that the three contributions assemble into
$\sum_{k=0}^{n-1}\frac{(-1)^kC_kC_{k+1}}{(k+1)!}\binom{n}{k+1}$.
\end{proof}

The first values of $F_n$ are
$$
1,\ 1,\ \tfrac53,\ \tfrac74,\ \tfrac{119}{60},\ \dots
$$
(for $n=1,2,3,4,5$); equivalently $n!\,F_n=1,2,10,42,238,1644,12528,\dots$.

\subsection{Boxes with a given hook length.}\label{sec:hookh}
By Lemma~\ref{lem:corner_hook_one}, $\co(\la)$ is the number of boxes of hook length $1$. In this subsection we extend Theorem~\ref{th:exactF2_sec1} to an arbitrary hook length $h\geq1$. Let
$$
\co_h(\la):=\#\{\square\in\la:\ h_\square=h\},\qquad
F_n^{(h)}:=\mathbb{E}[\co_h(\la_{(n)}^{(2)})]=\frac1{n!}\sum_{|\la|=n}(f_\la)^2\,\co_h(\la),
$$
so that $\co_1=\co$ and $F_n^{(1)}=F_n$. In the proof of Theorem~\ref{th:exactF2_sec1}, the product $\prod_{2\le j\le n}(h_\square^2-j^2)$ detects the boxes with $h_\square=1$, since it vanishes at the nodes $j^2$ for $j=2,\dots,n$. For a general $h$ we replace it by the Lagrange interpolation polynomial that vanishes at every node $j^2$ except $j=h$.

\begin{thm}\label{thm:exacth}
For all integers $h\ge1$ and $n\ge1$,
\begin{equation}\label{eq:exacth}
\begin{aligned}
F_n^{(h)}&=\sum_{m=h}^{n}(-1)^{m-h}\binom{n}{m}\,
\frac{h^2\,(m+1)!}{m\,(m-h)!\,(m+h)!}\,C_{m-1}C_m\\
&=2h^2\sum_{m=h}^{n}(-1)^{m-h}\binom{n}{m}\frac{\Lambda_{m-1}}{(m-h)!\,(m+h)!},
\end{aligned}
\end{equation}
where $\Lambda_{m-1}=\dfrac{(2m-2)!\,(2m-1)!}{(m-1)!\,m!^2}$ are the constants of Lemma~\ref{lem:okada}. For $h=1$ this is Theorem~\ref{th:exactF2_sec1}.
\end{thm}

For a polynomial $P$ let
$$
\Phi_P(n):=\mathbb{E}\Bigl[\sum_{\square\in\la}P(h_\square^2)\Bigr]
$$
be its Plancherel average over the boxes; it is linear in $P$. In the Newton basis
\begin{equation}\label{eq:newtonbasis}
N_r(x):=\prod_{j=1}^{r}(x-j^2)\qquad(r\ge0,\ N_0\equiv1)
\end{equation}
attached to the nodes $x=j^2$, Lemma~\ref{lem:okada} reads $\Phi_{N_r}(n)=A_r(n)=\Lambda_r\binom n{r+1}$ for $r\geq1$; this also holds for $r=0$, since $\Phi_{N_0}(n)=n$ and $\Lambda_0=1$.

\begin{lem}\label{lem:filter}
Fix $h\ge1$ and $n\ge1$, and let $Q_h$ be the unique polynomial of degree at most $n-1$ such that
\begin{equation}\label{eq:Qnodes}
Q_h(j^2)=\delta_{jh}\qquad(j=1,\dots,n).
\end{equation}
Then $\sum_{\square\in\la}Q_h(h_\square^2)=\co_h(\la)$ for every partition $\la$ of size $n$, and consequently $F_n^{(h)}=\Phi_{Q_h}(n)$.
\end{lem}

\begin{proof}
Every box of a partition of size $n$ satisfies $1\le h_\square\le n$, so $h_\square^2$ is one of the nodes $1^2,\dots,n^2$. By \eqref{eq:Qnodes}, $Q_h(h_\square^2)$ equals $1$ if $h_\square=h$ and $0$ otherwise. Summing over the boxes and taking the expectation gives the claim.
\end{proof}

In the Newton basis, the interpolation identity behind Lemma~\ref{lem:filter} appears as \cite[Eq.~(3.1)]{Schachinger2023}. There it is combined with Lemma~\ref{lem:okada} to compute the Plancherel average of $\sum_{\square}\log h_\square$.

\begin{lem}\label{lem:divdiff}
Write $Q_h=\sum_{r=0}^{n-1}c_r^{(h)}N_r$ in the basis \eqref{eq:newtonbasis}. Then
\begin{equation}\label{eq:cr}
c_r^{(h)}
=\begin{cases}
\displaystyle\Bigl(\prod_{\substack{1\le m\le r+1\\ m\ne h}}(h^2-m^2)\Bigr)^{-1}, & r\ge h-1,\\[2mm]
0,& r\le h-2.
\end{cases}
\end{equation}
In particular $c_r^{(h)}$ does not depend on $n$.
\end{lem}

\begin{proof}
Since $N_r$ is the Newton basis at the nodes $1^2,2^2,\dots$, the coefficient $c_r^{(h)}$ is the divided difference of $Q_h$ on the first $r+1$ nodes:
$$
c_r^{(h)}=Q_h[\,1^2,\dots,(r+1)^2\,]
=\sum_{i=1}^{r+1}\frac{Q_h(i^2)}{\prod_{1\le m\le r+1,\,m\ne i}(i^2-m^2)}
=\sum_{i=1}^{r+1}\frac{\delta_{ih}}{\prod_{1\le m\le r+1,\,m\ne i}(i^2-m^2)}.
$$
This is the displayed reciprocal product when $h\le r+1$, and the empty sum $0$ when $h>r+1$.
\end{proof}

\begin{lem}\label{lem:denom}
For integers $M\ge h\ge1$,
\begin{equation}\label{eq:denom}
\prod_{\substack{1\le m\le M\\ m\ne h}}(h^2-m^2)=\frac{(-1)^{M-h}}{2h^2}\,(M-h)!\,(M+h)!.
\end{equation}
\end{lem}

\begin{proof}
Write $h^2-m^2=(h-m)(h+m)$. First, $\prod_{m=1}^{h-1}(h-m)=(h-1)!$ and $\prod_{m=h+1}^{M}(h-m)=(-1)^{M-h}(M-h)!$. Second, $\prod_{m=1}^{M}(h+m)=(h+M)!/h!$, and removing the factor $m=h$, equal to $2h$, gives $\prod_{m\ne h}(h+m)=(h+M)!/(2h\cdot h!)$. Multiplying and using $(h-1)!/(2h\cdot h!)=1/(2h^2)$ yields \eqref{eq:denom}.
\end{proof}

\begin{proof}[Proof of Theorem~\ref{thm:exacth}]
By Lemmas~\ref{lem:filter} and~\ref{lem:divdiff} and the linearity of $\Phi$,
$$
F_n^{(h)}=\Phi_{Q_h}(n)=\sum_{r=0}^{n-1}c_r^{(h)}\,\Phi_{N_r}(n)
=\sum_{r=h-1}^{n-1}c_r^{(h)}\Lambda_r\binom n{r+1}.
$$
Put $m=r+1$, so that $m$ runs from $h$ to $n$. By \eqref{eq:cr} and Lemma~\ref{lem:denom} with $M=m$,
$$
c_{m-1}^{(h)}=\frac{2h^2\,(-1)^{m-h}}{(m-h)!\,(m+h)!},
$$
which gives the second expression in \eqref{eq:exacth}. For the first one, note that
$$
\frac{\Lambda_{m-1}}{C_{m-1}C_m}
=\frac{(2m-2)!\,(2m-1)!}{(m-1)!\,m!^2}\cdot\frac{(m-1)!\,m!\cdot m!\,(m+1)!}{(2m-2)!\,(2m)!}
=\frac{(m+1)!}{2m},
$$
so that $2h^2\Lambda_{m-1}=h^2(m+1)!\,C_{m-1}C_m/m$. For $h=1$ the factor $\frac{(m+1)!}{m\,(m-1)!\,(m+1)!}$ equals $\frac1{m!}$, and \eqref{eq:exacth} becomes Theorem~\ref{th:exactF2_sec1} with $m=k+1$.
\end{proof}

Since every box has exactly one hook length, $\sum_{h\geq1}F_n^{(h)}=\mathbb{E}[|\la|]=n$.

The simplest new case is $h=2$. Then $\frac{4(m+1)!}{m\,(m-2)!\,(m+2)!}=\frac{4(m^2-1)}{(m+2)!}$, so that
\begin{equation}\label{eq:Fn2}
F_n^{(2)}=\sum_{m=2}^{n}(-1)^{m}\binom{n}{m}\,\frac{4\,(m^2-1)}{(m+2)!}\,C_{m-1}C_m.
\end{equation}
The first values are $F_n^{(2)}=0,\,1,\,\tfrac13,\,\tfrac76,\,\tfrac{13}{10},\,\tfrac{73}{60},\,\tfrac{199}{140},\,\tfrac{1147}{672}$ for $n=1,\dots,8$.

\subsection{Linear recurrences.}\label{ssec:recs}
The summand of Theorem~\ref{th:exactF2_sec1} is a proper hypergeometric term, and Zeilberger's algorithm gives an inhomogeneous recurrence of order $4$ for $F_n$.

\begin{prop}\label{prop:rec}
The sequence $F_n$ satisfies, for all $n\geq 0$ (with $F_0=0$), the inhomogeneous linear recurrence
\begin{equation}\label{eq:rec}
\sum_{i=0}^{4} p_i(n)\,F_{n+i} = -2,
\end{equation}
with polynomial coefficients
$$
\begin{aligned}
p_0(n)&=-(n+1)(n+2)(n+3), &\qquad
p_1(n)&=(4n-5)(n+2)(n+3),\\
p_2(n)&=-(6n^2+7n-16)(n+3), &\qquad
p_3(n)&=4n^3+29n^2+73n+72,\\
p_4(n)&=-(n+4)^2(n+5). &&
\end{aligned}
$$
\end{prop}

\begin{proof}
The summand
$T(n,k):=(-1)^k\frac{C_kC_{k+1}}{(k+1)!}\binom{n}{k+1}$ of Theorem~\ref{th:exactF2_sec1} is a proper hypergeometric term, so Zeilberger's creative telescoping algorithm \cite{Petkovsek1996etal} produces a recurrence together with a certificate. Put
$$
Q_{n,k}:=\sum_{i=0}^4 p_i(n)\,T(n+i,k),\qquad
G_{n,k}:=\frac{(k+1)^3(k+2)}{n+4}\,T(n+4,k).
$$
We claim that $Q_{n,k}=G_{n,k+1}-G_{n,k}$ for all $n,k\geq 0$. If $k\geq n+4$, all the terms vanish. If $0\le k\le n+3$, then $T(n+4,k)\ne0$, and dividing by it, the claim becomes an identity between rational functions of $n$ and $k$, because
$$
\frac{T(n+i,k)}{T(n+4,k)}=\prod_{j=i+1}^{4}\frac{n+j-1-k}{n+j},\qquad
\frac{T(n+4,k+1)}{T(n+4,k)}=-\frac{4(2k+1)(2k+3)(n+3-k)}{(k+2)^3(k+3)};
$$
this rational identity is verified by a direct computation.
Summing $Q_{n,k}=G_{n,k+1}-G_{n,k}$ over $0\le k\le n+3$ telescopes to $G_{n,n+4}-G_{n,0}$. Here $G_{n,n+4}=0$ since $T(n+4,n+4)=0$, and $G_{n,0}=\frac{2}{n+4}\,T(n+4,0)=2$. On the other hand $T(n+i,k)=0$ for $k\ge n+i$, so that $\sum_{k=0}^{n+3}T(n+i,k)=F_{n+i}$ and $\sum_{k=0}^{n+3}Q_{n,k}=\sum_{i=0}^4p_i(n)F_{n+i}$. This proves \eqref{eq:rec}.
\end{proof}

The exact formula \eqref{eq:exacth} yields a linear recurrence for $F_n^{(h)}$ whose coefficients are polynomials in $n$ and $h$, so that $h$ is a parameter.

\begin{prop}\label{prop:rechall}
For every integer $h\geq1$ and every $n\geq0$ (with $F_0^{(h)}=0$),
\begin{equation}\label{eq:rechall}
\sum_{i=0}^{5}r_i(n)\,F_{n+i}^{(h)}=0,
\end{equation}
where
$$
\begin{aligned}
r_0(n)&=(n+1)(n+2)(n+3)(n+4),\\
r_1(n)&=-(n+2)(n+3)(n+4)(5n-1),\\
r_2(n)&=(n+3)(n+4)(10n^2+22n-19-h^2),\\
r_3(n)&=-(n+4)\bigl(10n^3+72n^2+(149-3h^2)n+72-12h^2\bigr),\\
r_4(n)&=5n^4+74n^3+(418-3h^2)n^2+(1081-27h^2)n+1093-61h^2,\\
r_5(n)&=-(n+5)^2(n+5-h)(n+5+h).
\end{aligned}
$$
\end{prop}

\begin{proof}
For $m\in\mathbb Z$ put $a_m:=2h^2(-1)^{m-h}\Lambda_{m-1}/\bigl((m-h)!\,(m+h)!\bigr)$ if $m\geq h$ and $a_m:=0$ if $m<h$, so that $F_n^{(h)}=\sum_{m\in\mathbb Z}\binom nm a_m$ for all $n\geq0$ by Theorem~\ref{thm:exacth}. From $\Lambda_m/\Lambda_{m-1}=4m(2m-1)(2m+1)/(m+1)^2$ we get
\begin{equation}\label{eq:amrec}
(m+1)^2\bigl((m+1)^2-h^2\bigr)\,a_{m+1}+4m(4m^2-1)\,a_m=0\qquad(m\in\mathbb Z):
\end{equation}
for $m\geq h$ this is the ratio $a_{m+1}/a_m$; for $m=h-1$ both terms vanish, since $(m+1)^2=h^2$ and $a_{h-1}=0$; for $m\le h-2$ we have $a_m=a_{m+1}=0$. Let $A(y)=\sum_m a_my^m$ and $\theta_y=y\,\frac{d}{dy}$. Multiplying \eqref{eq:amrec} by $y^{m+1}$ and summing over $m$ gives
$$
\theta_y^2\bigl(\theta_y^2-h^2\bigr)A+4y\,\theta_y\bigl(4\theta_y^2-1\bigr)A=0.
$$
On the other hand,
$$
\Phi(x):=\sum_{n\geq0}F_n^{(h)}x^n=\sum_m a_m\sum_{n}\binom nm x^n=\sum_m a_m\frac{x^m}{(1-x)^{m+1}}=\frac{1}{1-x}\,A\Bigl(\frac{x}{1-x}\Bigr).
$$
With $y=x/(1-x)$ we have $\theta_y=x(1-x)\frac{d}{dx}$. Substituting $A(y)=(1-x)\Phi(x)$ in the differential equation for $A$, multiplying by $1-x$ and removing the factor $x(x-1)$, a direct computation gives $\sum_{k=0}^4P_k(x)\,\Phi^{(k)}(x)=0$ with
$$
\begin{aligned}
P_4&=x^3(x-1)^5,\qquad P_3=2x^2(x-1)^3(8x^2-3x+3),\\
P_2&=x(x-1)^2\bigl(72x^3+18x^2+(13-h^2)x+h^2-7\bigr),\\
P_1&=(x-1)\bigl(96x^4+84x^3-(4h^2+56)x^2+(5h^2-17)x+1-h^2\bigr),\\
P_0&=24x^4+36x^3-(2h^2+46)x^2+3(h^2-1)x+1-h^2.
\end{aligned}
$$
Taking the coefficient of $x^{n+4}$, with $[x^N]\,x^a\Phi^{(b)}(x)=(N-a+1)(N-a+2)\cdots(N-a+b)\,F^{(h)}_{N-a+b}$, gives \eqref{eq:rechall}.
\end{proof}

\emph{The case $h=1$.} For $h=1$ every coefficient $r_i(n)$ is divisible by $n+4$, and after division \eqref{eq:rechall} becomes, for the expected number of distinct parts $F_n=F_n^{(1)}$,
\begin{equation}\label{eq:rech1}
\begin{aligned}
&(n+1)(n+2)(n+3)\,F_n-(n+2)(n+3)(5n-1)\,F_{n+1}+2(n+3)(5n^2+11n-10)\,F_{n+2}\\
&\quad-2(n+3)(5n^2+21n+10)\,F_{n+3}+(5n^3+54n^2+199n+258)\,F_{n+4}\\
&\quad-(n+5)^2(n+6)\,F_{n+5}=0.
\end{aligned}
\end{equation}
This is the difference of two consecutive instances of the inhomogeneous recurrence \eqref{eq:rec}: if $\mathcal L$ is the operator of Proposition~\ref{prop:rec} and $S$ is the shift $n\mapsto n+1$, the left-hand side of \eqref{eq:rech1} is $\bigl((S-1)\mathcal LF\bigr)(n)=(\mathcal LF)(n+1)-(\mathcal LF)(n)$, which vanishes because $\mathcal LF=-2$ is constant. In other words, for $h=1$ the operator of \eqref{eq:rechall} factors as $(n+4)(S-1)\circ\mathcal L$, and Proposition~\ref{prop:rechall} recovers Proposition~\ref{prop:rec} up to the additive constant, which is fixed by the initial values.

The leading terms of $r_0,\dots,r_5$ are $n^4,-5n^4,10n^4,-10n^4,5n^4,-n^4$, so that $\rho=1$ is the only characteristic root, with multiplicity $5$, for every $h$. Note that $r_5(n)$ vanishes for $n=h-5$, so that for $h\geq5$ the recurrence \eqref{eq:rechall} determines $F_{n+5}^{(h)}$ from the previous values only for $n\neq h-5$; the missing value is $F_h^{(h)}=2h^2\Lambda_{h-1}/(2h)!$.

\subsection{Integral representations.}\label{ssec:intrep}
The sum in Theorem~\ref{th:exactF2_sec1} is strongly alternating, with terms far larger than the total, so it does not directly reveal the growth of $F_n$. We therefore first rewrite $F_n$ as an integral, using the classical moment representation of the Catalan numbers.

\begin{lem}[{\cite[Eq.~(10)]{Penson2001Sixdeniers}}]\label{lem:moments}
For every integer $k\ge0$,
$$
C_k=\frac1{2\pi}\int_0^4 t^k\,w(t)\,dt,
\qquad
w(t)=\sqrt{\frac{4-t}{t}}.
$$
\end{lem}

Let
$$
L_n(u)=\sum_{m=0}^n\binom nm\frac{(-u)^m}{m!}
$$
be the Laguerre polynomial of degree $n$. We need two classical facts on $L_n$. The first one is the Mehler--Heine formula \cite[Theorem~8.1.3]{Szego1975}:
\begin{equation}\label{eq:MH}
\lim_{n\to\infty} L_n(u/n)=J_0(2\sqrt u)
\end{equation}
uniformly for $u$ in compact subsets of $[0,\infty)$, where $J_0(v)=\sum_{m\ge0}\frac{(-1)^m}{m!^2}\bigl(\frac v2\bigr)^{2m}$ is the Bessel function of order $0$, so that $J_0(2\sqrt u)=\sum_{m\ge0}\frac{(-u)^m}{m!^2}$. The second one is the inequality \cite[22.14.12--13]{Abramowitz1964Stegun}
\begin{equation}\label{eq:szego}
|L_n(u)|\le e^{u/2},\qquad |L^{(1)}_{n-1}(u)|\le n\,e^{u/2}\qquad (u\geq 0),
\end{equation}
where $L^{(1)}_{n-1}(u)=\sum_{m=0}^{n-1}\binom{n}{m+1}\frac{(-u)^m}{m!}$ is the generalized Laguerre polynomial; it satisfies $L_n'(u)=-L^{(1)}_{n-1}(u)$.

\begin{lem}\label{lem:rep}
For $n\geq 1$ we have
$$
F_n
=\frac1{2\pi}\int_0^4\frac{w(x)}{x}\bigl(1-P_n(x)\bigr)\,dx,
\qquad
P_n(x)=\sum_{m=0}^n\binom nm\frac{(-x)^m C_m}{m!}=\frac1{2\pi}\int_0^4 w(t)\,L_n(xt)\,dt.
$$
\end{lem}

\begin{proof}
By Lemma~\ref{lem:moments}, $C_{m-1}=\frac1{2\pi}\int_0^4 x^{m-1}w(x)\,dx$. Inserting this
into Theorem~\ref{th:exactF2_sec1}, written in the form
$F_n=\sum_{m=1}^n(-1)^{m-1}\binom nm\frac{C_{m-1}C_m}{m!}$, and interchanging the finite sum
with the integral, we get
$$
F_n
=\frac1{2\pi}\int_0^4\frac{w(x)}{x}
\sum_{m=1}^n(-1)^{m-1}\binom nm\frac{C_m}{m!}x^m\,dx
=\frac1{2\pi}\int_0^4\frac{w(x)}{x}\bigl(1-P_n(x)\bigr)\,dx,
$$
since the term $m=0$ of $P_n(x)$ is $C_0=1$. The integral converges, because $1-P_n(x)=\O(x)$ and $w(x)/x=\O(x^{-3/2})$ as $x\to 0$. The second expression of $P_n$ follows by writing $C_m=\frac1{2\pi}\int_0^4 t^m w(t)\,dt$.
\end{proof}

\section{The leading asymptotics}\label{sec:beta2}

In this section we deduce the leading term of $F_n^{(h)}$ from the exact formulas of Section~\ref{sec:exact}, and thereby prove part~(b) of Theorems~\ref{mth:full} and~\ref{mth:hookh} for $K_1$ and $K_1^{(h)}$. For $h=1$: in \eqref{eq:fullexp} we have $S_1(n)\sim\sqrt n$, while the other terms are $\O(1)$, so $K_1=\lim_{n\to\infty}F_n/\sqrt n$, and we must show that
\begin{equation}\label{eq:asyF2}
\mathbb{E}[\co(\la_{(n)}^{(2)})] \sim \frac{32}{3\pi^2}\sqrt{n}\qquad(n\to\infty).
\end{equation}
The limit \eqref{eq:asyF2} is due to Bufetov \cite{Bufetov2012}: it follows from his convergence in probability \cite[Lemma~2.1]{Bufetov2012} for $k=1$, because $0\le\co(\la)/\sqrt n<\sqrt2$ (Lemma~\ref{lem: upper_bound}). Bufetov's proof is local: he writes the number of corners as a sum of local patterns in the Maya encoding of $\la$, and uses the determinantal structure of the Poissonized Plancherel measure (the discrete Bessel process), its local convergence to the discrete sine process, a decay-of-correlation estimate and depoissonization. The constant $\frac{32}{3\pi^2}$ then appears as the integral of the two-point function of the sine process over the limit shape (see the Remark at the end of Section~\ref{sec:discussion}). His method gives the stronger convergence in probability, for every local pattern, but no lower-order terms. We give a different proof below, which starts from the exact formula of Theorem~\ref{th:exactF2_sec1}, uses no Poissonization, and needs only the Catalan moment sequence and classical properties of Laguerre polynomials; here $\frac{32}{3\pi^2}$ comes from a Bessel integral against the Catalan weight. It concerns only the expectation, but the same exact formula gives access to the lower-order terms (Sections~\ref{sec:rec} and~\ref{sec:constants}).

\subsection{Distinct parts.}\label{ssec:distinct}

\begin{proof}[Proof of Theorem~\ref{mth:full}\textup{(b)} for $K_1$]
We prove \eqref{eq:asyF2}; the argument uses only the exact formula of Theorem~\ref{th:exactF2_sec1}, through Lemma~\ref{lem:rep}.
Substitute $x=\xi/n$ in Lemma~\ref{lem:rep}. Since
$$
\frac{w(\xi/n)}{\xi/n}=2n^{3/2}\,\xi^{-3/2}\sqrt{1-\frac{\xi}{4n}}
\qquad\text{and}\qquad dx=\frac{d\xi}{n},
$$
we obtain
\begin{equation}\label{eq:ast}
\frac{F_n}{\sqrt n}
=\frac{1}{\pi}\int_0^{\infty}g_n(\xi)\,d\xi,\qquad
g_n(\xi)=\mathbf 1_{\{\xi\le 4n\}}\,\xi^{-3/2}\sqrt{1-\tfrac{\xi}{4n}}\,
\bigl(1-P_n(\xi/n)\bigr).
\end{equation}

\emph{Pointwise limit.} By Lemma~\ref{lem:rep}, $P_n(\xi/n)=\frac1{2\pi}\int_0^4 w(t)L_n(\xi t/n)\,dt$. For fixed $\xi\ge0$, \eqref{eq:MH} holds uniformly for $u=\xi t\in[0,4\xi]$, hence
$$
\lim_{n\to\infty}P_n(\xi/n)=
\tilde P(\xi):=\frac1{2\pi}\int_0^4 w(t)\,J_0\bigl(2\sqrt{\xi t}\bigr)\,dt
$$
and $\lim_{n\to\infty} g_n(\xi)=\xi^{-3/2}\bigl(1-\tilde P(\xi)\bigr)$.

\emph{Domination.} Let $0\le\xi\le 4n$ and $0\le t\le 4$, and put $u=\xi t/n\in[0,16]$. By \eqref{eq:szego},
$|1-L_n(u)|\le 1+e^{8}$ and
$|1-L_n(u)|=\bigl|\int_0^u L^{(1)}_{n-1}(s)\,ds\bigr|\le n u\, e^{8}=e^{8}\xi t$. Since $\frac1{2\pi}\int_0^4 w(t)\,dt=C_0=1$ and $\frac1{2\pi}\int_0^4 t\,w(t)\,dt=C_1=1$, we get
$|1-P_n(\xi/n)|\le\frac1{2\pi}\int_0^4 w(t)\,|1-L_n(\xi t/n)|\,dt\le 2e^{8}\min(\xi,1)$. Hence
$$
|g_n(\xi)|\le 2e^{8}\,\xi^{-3/2}\min(\xi,1),
$$
which is integrable on $(0,\infty)$.

\emph{The limit.} By dominated convergence, \eqref{eq:ast} yields
$$
\lim_{n\to\infty}\frac{F_n}{\sqrt n}=\frac1\pi\int_0^\infty\xi^{-3/2}\bigl(1-\tilde P(\xi)\bigr)\,d\xi
=\frac{1}{2\pi^2}\int_0^4 w(t)\int_0^\infty\xi^{-3/2}\bigl(1-J_0(2\sqrt{\xi t})\bigr)\,d\xi\,dt,
$$
where we used $\frac1{2\pi}\int_0^4w=1$ and Tonelli's theorem (the integrand is nonnegative since $|J_0|\le 1$). Integrating by parts and using $J_0'=-J_1$ and $\int_0^\infty J_1(v)\,v^{-1}\,dv=1$, we have
$\int_0^\infty v^{-2}\bigl(1-J_0(av)\bigr)\,dv=a$ for $a\geq 0$; with $\xi=v^2$ this gives
$$
\int_0^\infty\xi^{-3/2}\bigl(1-J_0(2\sqrt{\xi t})\bigr)\,d\xi=2\int_0^\infty v^{-2}\bigl(1-J_0(2\sqrt t\,v)\bigr)\,dv=4\sqrt t.
$$
Therefore
$$
\lim_{n\to\infty}\frac{F_n}{\sqrt n}=\frac{1}{2\pi^2}\int_0^4 w(t)\,4\sqrt t\,dt
=\frac{2}{\pi^2}\int_0^4 \sqrt{4-t}\,dt=\frac{2}{\pi^2}\cdot\frac{16}{3}=\frac{32}{3\pi^2}.
$$
This proves \eqref{eq:asyF2}, and hence part~(b) for $K_1$.
\end{proof}

\subsection{The leading asymptotics for every $h$.}\label{ssec:asyh}
The method of the previous subsection applies to the exact formula of Theorem~\ref{thm:exacth}. For $0\le s\le1$ and $n\geq1$ put
$$
G_n(s):=\sum_{m=1}^{n}(-1)^{m-1}\binom nm\frac{C_{m-1}C_m}{m!}\,s^m,
$$
so that $G_n(1)=F_n$ by Theorem~\ref{th:exactF2_sec1}.

\begin{lem}\label{lem:Gns}
There is a constant $C$ such that $|G_n(s)|\le C\sqrt{ns}$ for all $n\geq1$ and $0\le s\le1$. Moreover, for every $s\in[0,1]$,
$$
\lim_{n\to\infty}\frac{G_n(s)}{\sqrt n}=\frac{32}{3\pi^2}\sqrt s .
$$
\end{lem}

\begin{proof}
Both claims are clear for $s=0$; let $0<s\le1$. As in the proof of Lemma~\ref{lem:rep},
$G_n(s)=\frac1{2\pi}\int_0^4\frac{w(x)}{x}\bigl(1-P_n(sx)\bigr)\,dx$. Substitute $x=\xi/(ns)$. Since
$$
\frac{w(\xi/(ns))}{\xi/(ns)}=2(ns)^{3/2}\,\xi^{-3/2}\sqrt{1-\frac{\xi}{4ns}}
\qquad\text{and}\qquad dx=\frac{d\xi}{ns},
$$
we obtain
$$
\frac{G_n(s)}{\sqrt n}=\frac{\sqrt s}{\pi}\int_0^\infty g_{n,s}(\xi)\,d\xi,\qquad
g_{n,s}(\xi)=\mathbf 1_{\{\xi\le4ns\}}\,\xi^{-3/2}\sqrt{1-\tfrac{\xi}{4ns}}\,\bigl(1-P_n(\xi/n)\bigr).
$$
On the support of $g_{n,s}$ we have $\xi\le4ns\le4n$, so the domination step in the proof of \eqref{eq:asyF2} gives $|g_{n,s}(\xi)|\le2e^8\xi^{-3/2}\min(\xi,1)$. Hence $|G_n(s)|\le C\sqrt{ns}$ with $C=\frac{2e^8}{\pi}\int_0^\infty\xi^{-3/2}\min(\xi,1)\,d\xi=\frac{8e^8}{\pi}$. For fixed $s$, the pointwise limit of $g_{n,s}(\xi)$ is $\xi^{-3/2}\bigl(1-\tilde P(\xi)\bigr)$, as in the proof of \eqref{eq:asyF2}, and dominated convergence gives
$G_n(s)/\sqrt n\to\frac{\sqrt s}{\pi}\int_0^\infty\xi^{-3/2}\bigl(1-\tilde P(\xi)\bigr)\,d\xi=\frac{32}{3\pi^2}\sqrt s$.
\end{proof}

\begin{proof}[Proof of Theorem~\ref{mth:hookh}\textup{(b)} for $K_1^{(h)}$]
In \eqref{eq:fullexph} we have $S_1^{(h)}(n)\sim\sqrt n$, while the other terms are $\O(1)$. Hence $K_1^{(h)}=\lim_{n\to\infty}F_n^{(h)}/\sqrt n$, and we must show that, for every fixed $h\geq1$,
\begin{equation}\label{eq:asyh}
\mathbb{E}[\co_h(\la_{(n)}^{(2)})]\sim\frac{32h^2}{(4h^2-1)\pi^2}\sqrt n\qquad(n\to\infty).
\end{equation}
The proof of \eqref{eq:asyh} below uses only Theorem~\ref{thm:exacth}; it is used in the proof of part~(a) in Section~\ref{ssec:fullh}.

For $m\geq1$ let
$$
\rho_h(m):=h^2\,\frac{\prod_{j=1}^{h-1}(m-j)}{\prod_{j=2}^{h}(m+j)}.
$$
Since $\frac{m!}{(m-h)!}=m\prod_{j=1}^{h-1}(m-j)$ and $\frac{(m+1)!}{(m+h)!}=\prod_{j=2}^{h}(m+j)^{-1}$, we have $\frac{h^2(m+1)!}{m\,(m-h)!\,(m+h)!}=\frac{\rho_h(m)}{m!}$ for $m\ge h$, while $\rho_h(m)=0$ for $1\le m\le h-1$. Hence Theorem~\ref{thm:exacth} reads
$$
F_n^{(h)}=(-1)^{h-1}\sum_{m=1}^{n}(-1)^{m-1}\binom nm\frac{C_{m-1}C_m}{m!}\,\rho_h(m).
$$
The rational function $\rho_h(m)/h^2$ has numerator and denominator monic of the same degree $h-1$, and simple poles at $m=-2,\dots,-h$. Hence there are rational constants $\beta_2,\dots,\beta_h$ such that, for $m\geq1$,
\begin{equation}\label{eq:rhopf}
\rho_h(m)=h^2\Bigl(1+\sum_{j=2}^{h}\frac{\beta_j}{m+j}\Bigr)=h^2\Bigl(1+\sum_{j=2}^{h}\beta_j\int_0^1s^{m+j-1}\,ds\Bigr).
\end{equation}
Inserting this in the previous formula and exchanging the finite sum with the integrals, we get
\begin{equation}\label{eq:FhG}
F_n^{(h)}=(-1)^{h-1}h^2\Bigl(G_n(1)+\sum_{j=2}^h\beta_j\int_0^1s^{j-1}G_n(s)\,ds\Bigr).
\end{equation}
By Lemma~\ref{lem:Gns}, $|s^{j-1}G_n(s)|/\sqrt n\le C$ on $[0,1]$, and dominated convergence gives
$$
\begin{aligned}
\lim_{n\to\infty}\frac{F_n^{(h)}}{\sqrt n}
&=(-1)^{h-1}h^2\,\frac{32}{3\pi^2}\Bigl(1+\sum_{j=2}^h\beta_j\int_0^1s^{j-1/2}\,ds\Bigr)\\
&=(-1)^{h-1}h^2\,\frac{32}{3\pi^2}\Bigl(1+\sum_{j=2}^h\frac{\beta_j}{\frac12+j}\Bigr)
=(-1)^{h-1}\,\frac{32}{3\pi^2}\,\rho_h\bigl(\tfrac12\bigr),
\end{aligned}
$$
where $\rho_h(\frac12)$ is the value at $m=\frac12$ of the rational function $\rho_h$. Finally, since $\Gamma(h+\frac32)=(h+\frac12)(h-\frac12)\Gamma(h-\frac12)$ and $\Gamma(\frac52)=\frac34\Gamma(\frac12)$,
$$
\rho_h\bigl(\tfrac12\bigr)=h^2\,\frac{\prod_{j=1}^{h-1}(\frac12-j)}{\prod_{j=2}^{h}(\frac12+j)}
=(-1)^{h-1}h^2\,\frac{\Gamma(h-\frac12)\,\Gamma(\frac52)}{\Gamma(\frac12)\,\Gamma(h+\frac32)}
=(-1)^{h-1}\,\frac{3h^2}{4h^2-1},
$$
and the limit is $\frac{32}{3\pi^2}\cdot\frac{3h^2}{4h^2-1}=\frac{32h^2}{(4h^2-1)\pi^2}$. This proves \eqref{eq:asyh}, and hence part~(b) for $K_1^{(h)}$.
\end{proof}

For $h=1$, \eqref{eq:asyh} is \eqref{eq:asyF2}. The constants $\frac{32h^2}{(4h^2-1)\pi^2}$ decrease from $\frac{32}{3\pi^2}$ to $\frac8{\pi^2}$ as $h\to\infty$. The limit \eqref{eq:asyh} agrees with Bufetov's result \cite[Lemma~2.1]{Bufetov2012}.

\section{The complete asymptotic expansion}\label{sec:rec}

The limit \eqref{eq:asyF2} gives the leading term of $F_n=\mathbb{E}[\co(\la_{(n)}^{(2)})]$. In this section we prove part~(a) of Theorem~\ref{mth:full}. The exact formula of Theorem~\ref{th:exactF2_sec1} determines, through a linear recurrence, the \emph{shape} of the whole asymptotic expansion: the exponents, the logarithmic term and the additive constant are forced, and, once the leading constant $32/(3\pi^2)$ is known, only one connection constant remains free in the non-oscillatory part of the expansion (the oscillatory part involves a further complex constant). In particular, we shall obtain
\begin{equation}\label{eq:shape}
\mathbb{E}[\co(\la_{(n)}^{(2)})]=\frac{32}{3\pi^2}\Bigl(\sqrt n+\frac{3}{2^8}\,\frac{\log n}{\sqrt n}\Bigr)-\frac12+\frac{K_2}{\sqrt n}+\O\Bigl(\frac{\log n}{n^{3/2}}\Bigr).
\end{equation}

Note that $\frac{3}{2^8}\cdot\frac{32}{3\pi^2}=\frac{1}{8\pi^2}$, so the logarithmic term is $\frac{\log n}{8\pi^2\sqrt n}$. The recurrence of Proposition~\ref{prop:rec}, together with its initial values, determines the sequence $F_n$ and hence the constant $K_2$. The formal recurrence calculation does not itself evaluate this constant; an exact Mellin-transform argument in Section~\ref{ssec:k2evaluation} will give its value.

\subsection{Formal solutions.}
Let $\mathcal L$ be the operator $(\mathcal L y)(n)=\sum_{i=0}^4 p_i(n)\,y(n+i)$. Since $\sum_{i=0}^4p_i(n)=4$, the constant $-\frac12$ is a particular solution of \eqref{eq:rec}, and $G_n:=F_n+\frac12$ satisfies the homogeneous equation $\mathcal L G=0$. As $p_4(n)\neq0$ for $n\geq0$, the solutions of $\mathcal L G=0$ form a vector space of dimension~$4$.

The leading parts of the cubic coefficients $p_0,\dots,p_4$ are $-n^3,4n^3,-6n^3,4n^3,-n^3$, so the characteristic polynomial is $-(\rho-1)^4$: the only characteristic root is $\rho=1$, with multiplicity~$4$. To find the formal solutions we compute the action of $\mathcal L$ on powers. Expanding $(n+i)^\theta=n^\theta\sum_j\binom\theta j i^jn^{-j}$ one finds
\begin{equation}\label{eq:Lpow}
\mathcal L[n^\theta]=n^\theta\Bigl(c_0(\theta)+\frac{c_1(\theta)}{n}+\O(n^{-2})\Bigr),\qquad
c_0(\theta)=4-16\theta^2,\quad c_1(\theta)=-\theta(\theta^3+33\theta^2+24\theta-20),
\end{equation}
the terms of order $n^{\theta+3},n^{\theta+2},n^{\theta+1}$ cancelling because of the multiple root $\rho=1$.

\begin{prop}\label{prop:formal}
The equation $\mathcal LG=0$ has four formal solutions $S_1,S_2,Y,\bar Y$, where $S_1$ is a series in $n^{1/2-m}$ and $n^{1/2-m}\log n$ \textup($m\ge0$\textup) without a term in $n^{-1/2}$, $S_2$ is a series in $n^{-1/2-m}$ \textup($m\geq 0$\textup), $Y$ is $e^{8i\sqrt n}n^{-7/4}$ times a series in $n^{-1/2}$, and $\bar Y$ is the complex conjugate of $Y$. All coefficients are determined recursively, and the first ones are
\begin{align*}
S_1(n)&=n^{1/2}+\frac{3\log n}{2^{8}n^{1/2}}
+\frac{17329}{2^{18}n^{3/2}}-\frac{573\log n}{2^{17}n^{3/2}}
-\frac{2580007}{2^{26}n^{5/2}}+\frac{150693\log n}{2^{26}n^{5/2}}+\cdots,\\[3pt]
S_2(n)&=\frac1{n^{1/2}}-\frac{191}{2^{9}n^{3/2}}+\frac{50231}{2^{18}n^{5/2}}
-\frac{25906715}{2^{28}n^{7/2}}+\frac{11298377517}{2^{38}n^{9/2}}+\cdots,\\[3pt]
Y(n)&=e^{8i\sqrt n}n^{-7/4}
 \Bigl(1+\frac{365i}{192n^{1/2}}-\frac{146545}{73728n}
 -\frac{372015913i}{212336640n^{3/2}}
 +\frac{59893028257}{32614907904n^2}+\cdots\Bigr).
\end{align*}
\end{prop}

\begin{proof}
\emph{The algebraic solutions.} By \eqref{eq:Lpow}, the indicial equation is $c_0(\theta)=4-16\theta^2=0$, with roots $\theta=\pm\frac12$. For $\theta=-\frac12$, substituting $S_2=\sum_{m\ge0}a_mn^{-1/2-m}$ ($a_0=1$) into $\mathcal LS_2=0$ gives, at order $n^{-1/2-m}$, an equation $c_0(-\frac12-m)\,a_m=(\text{terms in }a_0,\dots,a_{m-1})$, which determines $a_m$ since $c_0(-\frac12-m)=-16m(m+1)\neq0$ for $m\geq1$.

For $\theta=\frac12$ the two indicial roots differ by the integer~$1$, and a logarithm appears. By \eqref{eq:Lpow}, $\mathcal L[n^{1/2}]=c_1(\tfrac12)\,n^{-1/2}+\O(n^{-3/2})=-\frac3{16}\,n^{-1/2}+\O(n^{-3/2})$. Differentiating \eqref{eq:Lpow} with respect to $\theta$ gives $\mathcal L[n^\theta\log n]=n^{\theta}\bigl(c_0(\theta)\log n+c_0'(\theta)\bigr)+\O(n^{\theta-1}\log n)$, hence, since $c_0(-\frac12)=0$ and $c_0'(-\frac12)=16$,
$$
\mathcal L\bigl[n^{-1/2}\log n\bigr]=16\,n^{-1/2}+\O(n^{-3/2}\log n).
$$
Therefore $\mathcal L\bigl[n^{1/2}+b\,n^{-1/2}\log n\bigr]=\O(n^{-3/2}\log n)$ exactly when $16b=\frac3{16}$, i.e., $b=\frac{3}{2^8}$. The higher coefficients of $n^{1/2-m}$ and $n^{1/2-m}\log n$ ($m\ge2$) are then determined recursively, because $c_0(\frac12-m)=16m(1-m)\neq0$ for $m\ge2$. (The coefficient of $n^{-1/2}$ in $S_1$ is free, as one may add a multiple of $S_2$; we take it equal to $0$.)

Explicitly, write $S(n)=\sum_{m\ge0}(a_m+b_m\log n)\,n^{\theta-m}$ and expand $\mathcal L S$ with
$$
(n+i)^{\theta-m}=n^{\theta-m}\sum_{j\ge0}\binom{\theta-m}{j}\frac{i^j}{n^{j}},\qquad
\log(n+i)=\log n+\sum_{j\ge1}\frac{(-1)^{j-1}}{j}\,\frac{i^j}{n^{j}}.
$$
The coefficients of $n^{\theta-m}\log n$ and $n^{\theta-m}$ give the triangular system
$$
c_0(\theta-m)\,b_m=\beta_m,\qquad c_0(\theta-m)\,a_m=\alpha_m-c_0'(\theta-m)\,b_m,
$$
where $\alpha_m,\beta_m$ are linear forms in $a_0,b_0,\dots,a_{m-1},b_{m-1}$ with rational coefficients. Solving it with $\theta=\frac12$, $a_0=1$, $b_0=a_1=0$, and with $\theta=-\frac12$, $a_0=1$, $b_0=0$ (then $b_m=0$ for all $m$), gives the displayed coefficients, which are rational numbers.

\emph{The oscillatory solutions.} The remaining two solutions are not algebraic. Substitute $y(n)=e^{\mu\sqrt n}n^{\kappa}$ and expand $\sqrt{n+i}-\sqrt n$ and $(1+i/n)^\kappa$ in powers of $n^{-1/2}$. A direct computation gives
$$
\frac{\mathcal L[y](n)}{y(n)}=-\frac{\mu^2(\mu^2+64)}{16}\,n-\frac{\mu\,(\mu^4+8\mu^2\kappa+72\mu^2+256\kappa+64)}{16}\,n^{1/2}+\O(1).
$$
The first coefficient vanishes for $\mu=\pm8i$ ($\mu=0$ corresponds to the algebraic solutions); then $\mu^4=4096$, $\mu^2=-64$, and the second coefficient vanishes if and only if $-448-256\kappa=0$, i.e., $\kappa=-\frac74$. The lower-order terms are then determined recursively as a series in $n^{-1/2}$; this gives the displayed coefficients of $Y$, and the root $\mu=-8i$ gives $\bar Y$, since $\mathcal L$ has real coefficients.
\end{proof}

\subsection{Proof of Theorem~\ref{mth:full}\textup{(a)}.}\label{ssec:proofshape}
We use the following classical result, in the formulation of \cite[Theorem~3.3]{Benoit2017JM}.

\begin{lem}[Birkhoff--Trjitzinsky \cite{Birkhoff1933Trjitzinsky}]\label{lem:BT}
Consider a linear recurrence $\sum_{i=0}^s b_i(n)\,y(n+i)=0$ with polynomial coefficients, $b_0b_s\ne0$. Call \emph{formal solution} a formal expression
$$
\bar y(n)=(n!)^{\kappa}\,\alpha^n\,e^{\pi(n)}\sum_{j=0}^{J}(\log n)^j\sum_{i\ge0}\beta_{j,i}\,n^{\theta-i/\rho}
$$
that satisfies the recurrence formally, where $\kappa\in\mathbb Q$, $\alpha\in\mathbb C\setminus\{0\}$, $\theta,\beta_{j,i}\in\mathbb C$, $\rho\geq1$ is an integer and $\pi(n)=\pi_1n^{1/\rho}+\dots+\pi_{\rho-1}n^{(\rho-1)/\rho}$. Then:
\begin{itemize}
\item[(i)] the recurrence has $s$ formally linearly independent formal solutions;
\item[(ii)] for any $s$ formally linearly independent formal solutions $\bar y_1,\dots,\bar y_s$, there exist solutions $y_1,\dots,y_s$ of the recurrence, defined for $n$ large enough and forming a basis of the space of solutions, such that $y_k(n)\sim\bar y_k(n)$ for every $k$; that is, for every $I\geq1$,
$$
y_k(n)=(n!)^{\kappa}\alpha^ne^{\pi(n)}\Bigl(\sum_{j=0}^{J}(\log n)^j\sum_{i<I}\beta_{j,i}\,n^{\theta-i/\rho}+\O\bigl(n^{\Re\theta-I/\rho}(\log n)^J\bigr)\Bigr).
$$
\end{itemize}
\end{lem}

The original proof of Birkhoff and Trjitzinsky \cite{Birkhoff1933Trjitzinsky} is long and was long considered incomplete; complete modern proofs are due to Immink \cite{Immink1991} (see also \cite{Put1997Singer}, and the discussion after \cite[Theorem~3.3]{Benoit2017JM}). See also \cite{Wimp1985Zeilberger} for applications to recurrences arising in combinatorics.

The four formal solutions $S_1,S_2,Y,\bar Y$ of Proposition~\ref{prop:formal} are of the above form (with $\kappa=0$, $\alpha=1$, and $\rho=1$ for $S_1,S_2$, $\rho=2$, $\pi(n)=\pm8i\sqrt n$ for $Y,\bar Y$). Their normalizations are fixed together: $S_1$ has leading term $n^{1/2}$ and no term in $n^{-1/2}$, $S_2$ has leading term $n^{-1/2}$, and $Y,\bar Y$ have leading terms $e^{\pm8i\sqrt n}n^{-7/4}$. They are formally linearly independent: in a formal relation $aS_1+bS_2+cY+d\bar Y=0$, the parts with the exponential factors $e^{8i\sqrt n}$, $e^{-8i\sqrt n}$ and $1$ vanish separately, which gives $c=d=0$, and then the coefficients of $n^{1/2}$ and $n^{-1/2}$ give $a=0$ and $b=0$. By Lemma~\ref{lem:BT}, each of them is the asymptotic expansion of an actual solution of $\mathcal LG=0$, and these four actual solutions form a basis of the solution space. Each actual solution is uniquely determined by its asymptotic expansion. Indeed, if a combination $aS_1+bS_2+cY+d\bar Y$ of the actual solutions is $\O(n^{-2})$, then $a=0$ and $b=0$ by the terms $n^{1/2}$ and $n^{-1/2}$, and $cY(n)+d\bar Y(n)=n^{-7/4}\bigl(ce^{8i\sqrt n}+de^{-8i\sqrt n}+\O(n^{-1/2})\bigr)$; the phases $8\sqrt n$ are dense modulo $2\pi$ (they are unbounded with increments tending to $0$), so the bracket does not tend to $0$ unless $c=d=0$. Hence the constants $K_1,K_2,K_3$ below are well defined. Since the recurrence has rational coefficients and the formal series $S_1,S_2$ have real coefficients, the actual solutions corresponding to $S_1,S_2$ may be taken real, and those corresponding to $Y,\bar Y$ complex conjugate. Hence $G_n=F_n+\frac12$ can be written as
$$
F_n+\frac12=K_1S_1(n)+K_2S_2(n)+K_3\,Y(n)+\bar K_3\,\bar Y(n)
$$
for real constants $K_1,K_2$ and a complex constant $K_3$, where $S_1,S_2,Y,\bar Y$ now denote the actual solutions. This is \eqref{eq:fullexp}, which proves part~(a). Since $|Y(n)|=\O(n^{-7/4})$ and $S_2(n)=\O(n^{-1/2})$, we have $F_n/\sqrt n\to K_1$, so that $K_1=\frac{32}{3\pi^2}$ by \eqref{eq:asyF2}. Inserting the expansions of Proposition~\ref{prop:formal} gives \eqref{eq:shape}, the oscillatory part $\O(n^{-7/4})$ being absorbed into the remainder. \qed

\subsection{The complete expansion.}
By part~(a), the four constants $K_1$, $K_2$, $K_3$ and~$\bar K_3$ determine $F_n$ completely, and $K_1$ and $K_2$ alone determine every non-oscillatory term of its asymptotic expansion. Indeed, applying Lemma~\ref{lem:BT}(ii) to $S_1$ and $S_2$ with $I=M+1$, and using $|Y(n)|=\O(n^{-7/4})$, we get from \eqref{eq:fullexp}, for every integer $M\geq1$,
\begin{equation}\label{eq:fullexpM}
\mathbb{E}[\co(\la_{(n)}^{(2)})]=-\frac12+\sum_{m=0}^{M}\frac{K_1\,(a_m+b_m\log n)+K_2\,a'_{m-1}}{n^{m-1/2}}+\O\bigl(n^{-7/4}\bigr)+\O\Bigl(\frac{\log n}{n^{M+1/2}}\Bigr),
\end{equation}
where $S_1(n)\sim\sum_m(a_m+b_m\log n)\,n^{1/2-m}$, $S_2(n)\sim\sum_m a'_m\,n^{-1/2-m}$ (with $a'_{-1}=0$) are the series of Proposition~\ref{prop:formal}.

Thus the oscillatory part has the additive expansion $K_3\,Y(n)+\bar K_3\,\bar Y(n)=2\Re\bigl(K_3\,e^{8i\sqrt n}\bigr)n^{-7/4}+\O(n^{-9/4})$. Its leading amplitude is $2|K_3|\,n^{-7/4}$ and its phase is $8\sqrt n+\arg K_3$. It first enters the expansion between the algebraic terms of order $n^{-3/2}$ and $n^{-5/2}$, so that \eqref{eq:fullexpM} with $M=2$ reads
$$
F_n=K_1\Bigl(\sqrt n+\frac{3\log n}{2^8\sqrt n}+\frac{17329}{2^{18}n^{3/2}}-\frac{573\log n}{2^{17}n^{3/2}}\Bigr)
+K_2\Bigl(\frac1{\sqrt n}-\frac{191}{2^9n^{3/2}}\Bigr)-\frac12+\O\bigl(n^{-7/4}\bigr).
$$
For instance, the coefficient of $n^{-3/2}\log n$ is $-\frac{573}{2^{17}}K_1=-\frac{191}{2^{12}\pi^2}$.

\subsection{The complete expansion for every $h$.}\label{ssec:fullh}
We now prove part~(a) of Theorem~\ref{mth:hookh}, which extends part~(a) of Theorem~\ref{mth:full} to every hook length. Throughout, we write $K_1^{(h)}=\frac{32h^2}{(4h^2-1)\pi^2}$, so that $K_1^{(1)}=K_1$; by \eqref{eq:asyh}, $K_1^{(h)}$ is the limit of $F_n^{(h)}/\sqrt n$. The constant term of the expansion is not determined by the recurrence \eqref{eq:rechall}, since the constants are solutions of it (indeed $\sum_{i=0}^5r_i(n)=0$). We obtain it from the case $h=1$, through the following observation: for $n\ge1$ and $0\le s\le1$,
\begin{equation}\label{eq:thinning}
G_n(s)=\sum_{k=0}^{n}\binom nk s^k(1-s)^{n-k}F_k .
\end{equation}
Indeed, write $F_k=\sum_m\binom km\alpha_m$ with $\alpha_m=(-1)^{m-1}C_{m-1}C_m/m!$. Since $\sum_k\binom nk\binom km s^k(1-s)^{n-k}=\binom nm s^m$, the right-hand side is equal to $\sum_m\binom nm\alpha_m s^m=G_n(s)$. In probabilistic terms, $G_n(s)=\mathbb{E}[F_X]$, where $X$ is a binomial random variable with parameters $n$ and~$s$.

\begin{lem}\label{lem:consth}
For every $h\geq1$,
$$
\lim_{n\to\infty}\Bigl(\mathbb{E}[\co_h(\la_{(n)}^{(2)})]-K_1^{(h)}\sqrt n\Bigr)=-\frac h2 .
$$
\end{lem}

\begin{proof}
By \eqref{eq:shape}, $F_k=K_1\sqrt k-\frac12+\varepsilon_k$ with $\varepsilon_k\to0$; let $\varepsilon^*=\sup_{k\ge0}|\varepsilon_k|<\infty$. Let $X$ be binomial with parameters $n$ and $s$, and $\mu=ns$. From the identity $\sqrt x-\sqrt\mu-\frac{x-\mu}{2\sqrt\mu}=-\frac{(\sqrt x-\sqrt\mu)^2}{2\sqrt\mu}\geq-\frac{(x-\mu)^2}{2\mu^{3/2}}$ for $x\ge0$, together with $\mathbb{E}[\sqrt X]\le\sqrt\mu$ and $\operatorname{Var}X\le\mu$, we get $|\mathbb{E}[\sqrt X]-\sqrt\mu|\le\frac1{2\sqrt\mu}$. By Chebyshev's inequality, $\mathbb{P}(X<\mu/2)\le4/\mu$. Hence, by \eqref{eq:thinning},
$$
\Bigl|G_n(s)-K_1\sqrt{ns}+\frac12\Bigr|\le\frac{K_1}{2\sqrt{ns}}+\sup_{k\ge ns/2}|\varepsilon_k|+\frac{4\varepsilon^*}{ns},
$$
which tends to $0$ uniformly for $s\in[n^{-1/2},1]$. For $s\in[0,n^{-1/2}]$ we use instead $|G_n(s)-K_1\sqrt{ns}|\le(C+K_1)\sqrt{ns}$, by Lemma~\ref{lem:Gns}.

By \eqref{eq:FhG} and the identity $K_1^{(h)}=(-1)^{h-1}h^2K_1\bigl(1+\sum_{j}\beta_j\int_0^1s^{j-1/2}\,ds\bigr)$ established in the proof of \eqref{eq:asyh},
$$
F_n^{(h)}-K_1^{(h)}\sqrt n=(-1)^{h-1}h^2\Bigl(\bigl(G_n(1)-K_1\sqrt n\bigr)+\sum_{j=2}^h\beta_j\int_0^1s^{j-1}\bigl(G_n(s)-K_1\sqrt{ns}\bigr)\,ds\Bigr).
$$
For $j\geq2$, the part of the integral over $[0,n^{-1/2}]$ is at most $(C+K_1)\sqrt n\int_0^{n^{-1/2}}s^{j-1/2}\,ds\le(C+K_1)\,n^{1/4-j/2}\to0$, and the part over $[n^{-1/2},1]$ tends to $-\frac12\int_0^1s^{j-1}ds=-\frac1{2j}$. Therefore
$$
\lim_{n\to\infty}\bigl(F_n^{(h)}-K_1^{(h)}\sqrt n\bigr)=-\frac12(-1)^{h-1}h^2\Bigl(1+\sum_{j=2}^h\frac{\beta_j}{j}\Bigr)=-\frac12(-1)^{h-1}\rho_h(0)
$$
by \eqref{eq:rhopf} at $m=0$. Finally $\rho_h(0)=h^2\prod_{j=1}^{h-1}(-j)\big/\prod_{j=2}^hj=(-1)^{h-1}h$.
\end{proof}

We now determine the formal solutions of \eqref{eq:rechall}, as in Proposition~\ref{prop:formal}. Let $\mathcal R$ be the operator $(\mathcal Ry)(n)=\sum_{i=0}^5r_i(n)\,y(n+i)$. Expanding $(n+i)^\theta$ in powers of $n^{-1}$ as in \eqref{eq:Lpow}, one finds
\begin{equation}\label{eq:Rpow}
\mathcal R[n^\theta]=n^\theta\Bigl(d_0(\theta)+\frac{d_1(\theta)}n+\O(n^{-2})\Bigr),
\end{equation}
where
$$
d_0(\theta)=-4\theta(4\theta^2-1),\qquad
d_1(\theta)=\theta\bigl(h^2\theta^2-\theta^4-40\theta^3-48\theta^2+46\theta-6\bigr).
$$

\begin{prop}\label{prop:formalh}
For every $h\geq1$, the equation $\mathcal Ry=0$ has the five formal solutions $1$,
$$
\begin{aligned}
S_1^{(h)}(n)&=\sqrt n+\frac{4h^2-1}{2^8}\,\frac{\log n}{\sqrt n}+\frac{A_h}{n^{3/2}}+\frac{B_h\log n}{n^{3/2}}+\O\Bigl(\frac{\log n}{n^{5/2}}\Bigr),\\
S_2^{(h)}(n)&=\frac1{\sqrt n}+\frac{4h^2-577}{3\cdot2^9}\,\frac1{n^{3/2}}+\O\Bigl(\frac1{n^{5/2}}\Bigr),\\
Y^{(h)}(n)&=e^{8i\sqrt n}\,n^{-7/4}\Bigl(1+\frac{(48h^2+317)\,i}{192\sqrt n}+\O\bigl(n^{-1}\bigr)\Bigr),
\end{aligned}
$$
where
$$
A_h=\frac{-368h^4+36376h^2+119953}{9\cdot2^{18}},\qquad
B_h=\frac{(4h^2-1)(4h^2-577)}{3\cdot2^{17}},
$$
$S_1^{(h)}$ is a series in $n^{1/2-m}$ and $n^{1/2-m}\log n$ $(m\ge0)$ without term in $n^{-1/2}$, $S_2^{(h)}$ is a series in $n^{-1/2-m}$ $(m\ge0)$, and $Y^{(h)}$ is $e^{8i\sqrt n}n^{-7/4}$ times a series in $n^{-1/2}$; the fifth formal solution is the complex conjugate $\overline{Y^{(h)}}$. The coefficients of $S_1^{(h)}$ and $S_2^{(h)}$ are polynomials in $h^2$ with rational coefficients. For $h=1$, $S_1^{(1)}$, $S_2^{(1)}$ and $Y^{(1)}$ are the formal solutions of Proposition~\ref{prop:formal}.
\end{prop}

\begin{proof}
The constant $1$ is a solution because $\sum_ir_i(n)=0$. The indicial polynomial $d_0(\theta)$ has the roots $0$ and $\pm\frac12$. For the root $\theta=0$ the solution is the constant $1$; for $\theta=-\frac12$ the coefficients of $S_2^{(h)}$ are determined recursively because $d_0(-\frac12-m)\neq0$ for $m\geq1$, exactly as in the proof of Proposition~\ref{prop:formal}. For $\theta=\frac12$, we have $\mathcal R[n^{1/2}]=d_1(\frac12)\,n^{-1/2}+\O(n^{-3/2})$ with $d_1(\frac12)=\frac{4h^2-1}{32}$, and $\mathcal R[n^{-1/2}\log n]=d_0'(-\frac12)\,n^{-1/2}+\O(n^{-3/2}\log n)$ with $d_0'(-\frac12)=-8$. Hence the coefficient of $n^{-1/2}\log n$ in $S_1^{(h)}$ is $\frac{4h^2-1}{256}$, and the next coefficients are determined recursively since $d_0(\frac12-m)\neq0$ for $m\geq2$ (the triangular system of the proof of Proposition~\ref{prop:formal}, with $d_0,d_1,\dots$ in place of $c_0,c_1,\dots$); solving it gives $A_h$, $B_h$ and the coefficient of $n^{-3/2}$ in $S_2^{(h)}$. For the oscillatory solutions, substituting $y(n)=e^{\mu\sqrt n}n^\kappa$ gives
$$
\frac{\mathcal R[y](n)}{y(n)}=-\frac{\mu^3(\mu^2+64)}{32}\,n^{3/2}-\frac{\mu^2\bigl(5\mu^4+40\kappa\mu^2+360\mu^2+1536\kappa+768\bigr)}{128}\,n+\O(n^{1/2}),
$$
so that $\mu=\pm8i$ and, with $\mu^2=-64$, the second coefficient vanishes if and only if $-1792-1024\kappa=0$, i.e., $\kappa=-\frac74$; neither value depends on $h$. Substituting $y(n)=e^{8i\sqrt n}n^{-7/4}(1+d_1n^{-1/2}+d_2n^{-1})$, the first nonzero coefficient of $\mathcal R[y](n)/y(n)$, at order $n^{1/2}$, is $-64i\bigl(h^2+4i\,d_1+\frac{317}{48}\bigr)$; hence $d_1=\frac{(48h^2+317)\,i}{192}$, and the root $\mu=-8i$ gives the complex conjugate $\overline{Y^{(h)}}$. The higher coefficients are determined recursively in the same way. Unlike the exponents, the coefficient $d_1$ depends on $h$, so that $Y^{(h)}\neq Y^{(1)}$ for $h\geq2$. Finally, for $h=1$ the operator $\mathcal R$ is $(n+4)(S-1)\circ\mathcal L$, so the formal solutions of $\mathcal LG=0$ are formal solutions of $\mathcal Ry=0$; with the same normalizations (leading coefficient $1$, and no term in $n^{-1/2}$ for $S_1$), they coincide.
\end{proof}

\begin{proof}[Proof of Theorem~\ref{mth:hookh}\textup{(a)}]
The five formal solutions of Proposition~\ref{prop:formalh} are normalized together by their leading terms $1$, $n^{1/2}$, $n^{-1/2}$ and $e^{\pm8i\sqrt n}n^{-7/4}$, with no term in $n^{-1/2}$ in $S_1^{(h)}$ (and no constant term in $S_1^{(h)}$ and $S_2^{(h)}$, which are series in $n^{1/2-m}$ and $n^{-1/2-m}$). They are formally linearly independent: as in Section~\ref{ssec:proofshape}, the parts with the factors $e^{8i\sqrt n}$, $e^{-8i\sqrt n}$ and $1$ of a formal relation vanish separately, and the non-oscillatory part then vanishes because its terms in $n^{1/2}$, $1$ and $n^{-1/2}$ come only from $S_1^{(h)}$, the constant $1$ and $S_2^{(h)}$, respectively. By Lemma~\ref{lem:BT}, they are the asymptotic expansions of five actual solutions of $\mathcal Ry=0$, forming a basis of the solution space. We may take the constant $1$ itself, real solutions $S_1^{(h)}$ and $S_2^{(h)}$, and complex conjugate solutions $Y^{(h)}$ and $\overline{Y^{(h)}}$. Since $F_n^{(h)}$ is a real solution of $\mathcal Ry=0$ (Proposition~\ref{prop:rechall}), we get
$$
F_n^{(h)}=K\,S_1^{(h)}(n)+K'\,S_2^{(h)}(n)+D+K_3^{(h)}\,Y^{(h)}(n)+\bar K_3^{(h)}\,\overline{Y^{(h)}}(n)
$$
for real constants $K,K',D$ and a complex constant $K_3^{(h)}$. Dividing by $\sqrt n$ and letting $n\to\infty$ gives $K=K_1^{(h)}$ by \eqref{eq:asyh}. Then $F_n^{(h)}-K_1^{(h)}\sqrt n\to D$, since $S_1^{(h)}(n)-\sqrt n$, $S_2^{(h)}(n)$ and $Y^{(h)}(n)$ tend to $0$; hence $D=-\frac h2$ by Lemma~\ref{lem:consth}. With $K_2^{(h)}=K'$, this is \eqref{eq:fullexph}. For $h=1$, the decomposition \eqref{eq:fullexph} coincides with \eqref{eq:fullexp}. Indeed, the actual solutions of $\mathcal LG=0$ used there are solutions of $\mathcal Ry=0$ with the same asymptotic expansions, and an actual solution of $\mathcal Ry=0$ is determined by its asymptotic expansion: a nonzero combination of $1,S_1^{(h)},S_2^{(h)},Y^{(h)},\overline{Y^{(h)}}$ is not $\O(n^{-2})$, by the argument of Section~\ref{ssec:proofshape} (the terms $n^{1/2}$, $1$, $n^{-1/2}$, and then the density of the phases $8\sqrt n$ modulo $2\pi$). In particular the constants $K_2^{(h)}$ and $K_3^{(h)}$ are well defined.
\end{proof}

By Proposition~\ref{prop:formalh}, the coefficient of $\frac{\log n}{\sqrt n}$ in $F_n^{(h)}$ is $K_1^{(h)}\frac{4h^2-1}{2^8}=\frac{h^2}{8\pi^2}$, and that of $\frac{\log n}{n^{3/2}}$ is $K_1^{(h)}B_h=\frac{h^2(4h^2-577)}{3\cdot2^{12}\pi^2}$; the coefficient of $n^{-3/2}$ is $K_1^{(h)}A_h+\frac{4h^2-577}{3\cdot2^9}K_2^{(h)}$. With part~(b), this gives \eqref{eq:Fhexplicit}. For $h=1$ we recover $A_1=\frac{17329}{2^{18}}$, $B_1=-\frac{573}{2^{17}}$ and $\frac{4-577}{3\cdot2^9}=-\frac{191}{2^9}$, as in Section~\ref{sec:rec}.

\section{The connection constants}\label{sec:constants}

The expansions of Section~\ref{sec:rec} involve the connection constants $K_1$, $K_2$ and $K_3$, which are determined by the recurrence and the initial values but are not evaluated by the formal calculation. The constant $K_1$ was obtained in Section~\ref{sec:beta2}. In this section we evaluate $K_2$ by a Mellin transform of the ordinary generating function, and $K_3$ by an exact Poissonization and a localized Fourier transform. Both calculations use the existence of the asymptotic decomposition proved in Section~\ref{sec:rec}. In Section~\ref{ssec:consth} we deduce the constants $K_2^{(h)}$ and $K_3^{(h)}$ for every hook length from $K_1$, $K_2$ and $K_3$.

\subsection{Evaluation of the second connection constant}
\label{ssec:k2evaluation}

We prove that, with the normalization of Proposition~\ref{prop:formal}, the constant $K_2$ in \eqref{eq:fullexp} and \eqref{eq:shape} is $\frac{4}{3\pi^2}+\frac{2\gamma+12\log 2-1}{8\pi^2}$.

\begin{proof}[Proof of Theorem~\ref{mth:full}\textup{(b)} for $K_2$]
We determine $K_2$ from the singular behavior of the ordinary generating
function.  Only an Abelian comparison is needed, since the existence of
the coefficient $K_2$ has already been proved in part~(a); see \eqref{eq:shape}.

\emph{A positive Laplace representation of the generating function.}
Define, for $0<u<16$,
\[
 \mu(u)=\frac{1}{4\pi^2}
       \int_{u/4}^{4} \frac{w(x)w(u/x)}{x^2}\,dx.
\]
Lemma~\ref{lem:rep}, followed by the substitution $u=xt$, gives
\begin{equation}\label{eq:k2mu}
 F_n=\int_0^{16}\mu(u)\bigl(1-L_n(u)\bigr)\,du.
\end{equation}
This is an absolutely convergent integral.  Indeed, $\mu$ is positive,
and, for every complex $r$ with $\Re r>1/2$, its moments are
\begin{equation}\label{eq:k2moments}
 \begin{split}
 M(r):=\int_0^{16}u^r\mu(u)\,du
 &=\left(\frac1{2\pi}\int_0^4x^{r-1}w(x)\,dx\right)
   \left(\frac1{2\pi}\int_0^4t^r w(t)\,dt\right)\\
 &=\frac{16^r}{4\pi}
   \frac{\Gamma(r-1/2)\Gamma(r+1/2)}
        {\Gamma(r+1)\Gamma(r+2)}.
 \end{split}
\end{equation}
Here the beta integral
\[
 \frac1{2\pi}\int_0^4x^v w(x)\,dx
   =\frac{4^v\Gamma(v+1/2)}{\sqrt\pi\,\Gamma(v+2)}
 \qquad (\Re v>-1/2)
\]
extends the Catalan moment formula.  In particular $M(1)=1$.

For $0<z<1$, put $h=1-z$, $q=z/(1-z)$, and
\[
 B(z)=\sum_{n\ge0}F_nz^n,
 \qquad H(q)=\int_0^{16}\mu(u)(1-e^{-qu})\,du.
\]
The Laguerre generating function gives the exact identity
\begin{equation}\label{eq:k2ogf}
 B(z)=\frac{H(q)}{1-z}.
\end{equation}
To justify exchanging the series with the integral despite the
nonintegrability of $\mu$ at zero, use
\[
 |1-L_n(u)|\le n e^8 u\qquad(0\le u\le16),
\]
which follows from $L_n'=-L_{n-1}^{(1)}$ and
\eqref{eq:szego}.  Thus
\[
 \int_0^{16}\mu(u)\sum_{n\ge0}z^n|1-L_n(u)|\,du
 \le e^8 M(1)\sum_{n\ge0}nz^n<\infty.
\]
The identity
$\sum_{n\ge0}L_n(u)z^n=(1-z)^{-1}\exp(-uz/(1-z))$
then proves \eqref{eq:k2ogf}.

\emph{Mellin asymptotics of $H$.}
For $-1<\Re s<-1/2$, Tonelli's theorem and then absolute
convergence for complex $s$ give
\[
 \int_0^\infty q^{s-1}H(q)\,dq
   =-\Gamma(s)M(-s)=:\Phi(s),
\]
because
$\int_0^\infty q^{s-1}(1-e^{-qu})\,dq=-\Gamma(s)u^{-s}$
for $-1<\Re s<0$.  Hence
\begin{equation}\label{eq:k2mellin}
 \Phi(s)=-\frac{16^{-s}}{4\pi}
 \frac{\Gamma(s)\Gamma(-s-1/2)\Gamma(1/2-s)}
      {\Gamma(1-s)\Gamma(2-s)}.
\end{equation}
All interchanges are justified, since $M(-\Re s)<\infty$ by
\eqref{eq:k2moments}.  Mellin inversion therefore yields
\[
 H(q)=\frac{1}{2\pi i}
       \int_{-3/4-i\infty}^{-3/4+i\infty}\Phi(s)q^{-s}\,ds.
\]
For completeness, the transform and its inverse are absolutely
integrable on this line: the preceding calculation gives integrability
of $q^{-7/4}H(q)$, and Stirling's formula gives, uniformly in every
closed vertical strip avoiding the finitely many poles,
\begin{equation}\label{eq:k2stirling}
 \Phi(\sigma+i\tau)
 =O\!\left((1+|\tau|)^{\sigma-7/2}
             e^{-\pi|\tau|/2}\right)
 \qquad (|\tau|\to\infty).
\end{equation}

Shift the line of integration to $\Re s=1$.  The horizontal integrals
vanish by \eqref{eq:k2stirling}; the new vertical integral is $O(q^{-1})$
for real $q\to+\infty$.  The poles crossed are $s=-1/2$, $0$ and
$1/2$.  Put
\[
 A=\frac{16}{3\pi^{3/2}},\qquad
 b=\frac1{8\pi^{3/2}},\qquad
 D=\gamma+10\log2-1.
\]
Their Laurent expansions are
\begin{align*}
 \Phi(s)&=-\frac{A}{s+1/2}+O(1) && (s\to-1/2),\\
 \Phi(s)&=\frac1{2s}+O(1) && (s\to0),\\
 \Phi(s)&=\frac{b}{(s-1/2)^2}
           -\frac{bD}{s-1/2}+O(1) && (s\to1/2).
\end{align*}
For the last identity one may use
\[
 \Gamma(-\varepsilon)=-\varepsilon^{-1}-\gamma+O(\varepsilon),
 \qquad
 \Gamma(-1-\varepsilon)=\varepsilon^{-1}+\gamma-1
                         +O(\varepsilon),
\]
together with
$\psi(1/2)=-\gamma-2\log2$ and
$\psi(3/2)=\psi(1/2)+2$.
Shifting to the right subtracts the residues of $\Phi(s)q^{-s}$.
Consequently,
\begin{equation}\label{eq:k2Hasymp}
 H(q)=A\sqrt q-\frac12
       +\frac{b(\log q+D)}{\sqrt q}+O(q^{-1}).
\end{equation}
This includes a proved remainder; no assertion about oscillatory
asymptotics or depoissonization has been used.

Substitute $q=(1-h)/h$ in \eqref{eq:k2ogf} and
\eqref{eq:k2Hasymp}.  As $h\to0^+$, this gives
\begin{equation}\label{eq:k2Bfromintegral}
 B(1-h)=A h^{-3/2}-\frac12h^{-1}
       +h^{-1/2}\left[-\frac A2+b(-\log h+D)\right]+O(1).
\end{equation}

\emph{Comparison with the already proved coefficient expansion.}
By \eqref{eq:shape},
\[
 F_n=K_1 n^{1/2}-\frac12
       +a\,n^{-1/2}\log n+K_2n^{-1/2}+R_n,
 \qquad
 a=\frac1{8\pi^2},\quad
 R_n=O(n^{-3/2}\log n).
\]
The remainder is absolutely summable.  With $t=-\log z$, the elementary
Abelian estimates
\begin{align*}
 \sum_{n\ge1}n^{1/2}e^{-tn}
   &=\Gamma(3/2)t^{-3/2}+O(1),\\
 \sum_{n\ge1}n^{-1/2}e^{-tn}
   &=\sqrt\pi\,t^{-1/2}+O(1),\\
 \sum_{n\ge1}n^{-1/2}\log n\,e^{-tn}
   &=\sqrt\pi\,t^{-1/2}\bigl(\psi(1/2)-\log t\bigr)+O(1)
\end{align*}
hold as $t\to0^+$.  They follow, for example, by comparing each sum
with its integral using the Euler--Maclaurin formula on $[1,\infty)$.
The error is uniformly bounded: the endpoint values and first
derivatives are bounded, the second derivatives have uniformly bounded
$L^1[1,\infty)$ norms, and the integrals over $(0,1)$ are bounded.
The three integrals over $(0,\infty)$ are exactly the displayed gamma
expressions.

Since $t=h+h^2/2+O(h^3)$, and since
$K_1\Gamma(3/2)=A$ and $a\sqrt\pi=b$, these estimates imply
\begin{equation}\label{eq:k2BfromBT}
 B(1-h)=A h^{-3/2}-\frac12h^{-1}
 +h^{-1/2}\left[-\frac{3A}{4}
             +b(-\log h+\psi(1/2))+K_2\sqrt\pi\right]+O(1).
\end{equation}
Comparing the coefficients of $h^{-1/2}$ in
\eqref{eq:k2Bfromintegral} and \eqref{eq:k2BfromBT} gives
\[
 K_2\sqrt\pi=\frac A4+b\bigl(D-\psi(1/2)\bigr).
\]
Using the values of $A,b,D$ and $\psi(1/2)$ proves
\[
 K_2=\frac{4}{3\pi^2}
      +\frac{2\gamma+12\log2-1}{8\pi^2}.
\]
\end{proof}

Consequently, after the constant term $-\frac12$ is separated, every coefficient of a power of $n$ and $\log n$ in the non-oscillatory expansion \eqref{eq:fullexpM} belongs to $\mathbb Q\pi^{-2}+\mathbb Q\gamma\pi^{-2}+\mathbb Q\pi^{-2}\log2$, because the formal coefficients in Proposition~\ref{prop:formal} are rational.

\subsection{The oscillatory constant $K_3$.}\label{ssec:cconst}
The normalization of $Y$ in Proposition~\ref{prop:formal}, with leading coefficient~$1$, fixes the constant
$K_3$ in \eqref{eq:fullexp}. We prove that
\[
 K_3=-\frac{e^{8+i\pi/4}}{1024\pi^{3/2}}.
\]
Consequently the oscillatory part has the additive expansion
\[
 K_3Y(n)+\bar K_3\bar Y(n)
 =-\frac{e^8}{512\pi^{3/2}}\,
       n^{-7/4}\cos\left(8\sqrt n+\frac\pi4\right)
       +\O(n^{-9/4}).
\]

\begin{proof}[Proof of Theorem~\ref{mth:full}\textup{(b)} for $K_3$]
We first Poissonize the exact formula, and then use a smooth Fourier
multiplier to isolate the upper endpoint.  This avoids the need for a
uniform remainder estimate in the Laguerre asymptotics.

\emph{The exact Poissonized integral and its endpoint density.}
For $t\geq0$ define
\[
 \mathcal P(t)=e^{-t}\sum_{n\geq0}F_n\frac{t^n}{n!}.
\]
Theorem~\ref{th:exactF2_sec1} and
$e^{-t}\sum_{n\geq m}\binom nm t^n/n!=t^m/m!$ give
\begin{equation}\label{eq:c-poisson}
 \mathcal P(t)=\sum_{m\geq1}
   \frac{(-1)^{m-1}C_{m-1}C_m}{(m!)^2}\,t^m
 =\int_0^{16}\mu(u)\bigl(1-J_0(2\sqrt{tu})\bigr)\,du,
\end{equation}
where $\mu$ is the positive density of \eqref{eq:k2mu}. A direct rewrite
of its defining integral gives
\begin{equation}\label{eq:c-mudensity}
 \mu(u)=\frac{1}{4\pi^2\sqrt u}
        \int_{u/4}^4
        \frac{\sqrt{(4-x)(4x-u)}}{x^{5/2}}\,dx,
       \qquad 0<u<16.
\end{equation}
By \eqref{eq:k2moments}, for every integer $m\geq1$ its moments are
$\int_0^{16}u^m\mu(u)\,du=C_{m-1}C_m$.
All the series interchanges here are absolute: the Catalan bound
$C_m\leq4^m$ bounds the resulting series by a constant times
$\sum_{m\geq1}(16t)^m/(m!)^2$.  In particular,
\begin{equation}\label{eq:c-firstmoment}
 \int_0^{16}u\mu(u)\,du=1.
\end{equation}

Put $s=16-u$ in \eqref{eq:c-mudensity} and set $x=4-sz/4$.
This gives the more precise endpoint formula
\[
 \mu(16-s)=\frac{s^2}{32\pi^2\sqrt{16-s}}
 \int_0^1\frac{\sqrt{z(1-z)}}{(4-sz/4)^{5/2}}\,dz
 =s^2 A(s),\qquad A(0)=\frac{1}{2^{15}\pi},
\]
where $A$ is real analytic near $0$.
Set
\[
 f(v)=\mathcal P(v^2/4),\qquad h(r)=2r\mu(r^2),\quad 0<r<4.
\]
Then for all real $v$,
\begin{equation}\label{eq:c-hankel}
 f(v)=\int_0^4h(r)\bigl(1-J_0(vr)\bigr)\,dr,
 \qquad \int_0^4r^2h(r)\,dr=1,
\end{equation}
and, near $r=4$,
\begin{equation}\label{eq:c-hendpoint}
 h(r)=(4-r)^2 H(4-r),\qquad H(0)=\frac{1}{64\pi},
\end{equation}
with $H$ real analytic near $0$.

\emph{Exact isolation of the positive endpoint frequency.}
Choose a real function $\phi\in C_c^\infty(\mathbb R)$ supported in
$(3,5)$ and equal to $1$ on a neighborhood of $4$, and put
\[
 \mathcal K(w)=\frac{1}{2\pi}\int_{\mathbb R}\phi(x)e^{ixw}\,dx,
 \qquad (\mathcal K*f)(v)=\int_{\mathbb R}\mathcal K(w)f(v-w)\,dw.
\]
Thus $\mathcal K$ is a Schwartz function and our Fourier convention is
$\int_{\mathbb R}\mathcal K(w)e^{-ixw}\,dw=\phi(x)$.
In particular,
\begin{equation}\label{eq:c-filtermoments}
 \int_{\mathbb R}w^j\mathcal K(w)\,dw=0\quad(j\geq0),\qquad
 \phi(4)=1,\quad\phi(-4)=0.
\end{equation}
The Bessel integral \cite[Eq.~9.1.18]{Abramowitz1964Stegun} reads
\[
 J_0(vr)=\frac1\pi\int_{-r}^r
                 \frac{e^{ivx}}{\sqrt{r^2-x^2}}\,dx.
\]
Apply it inside \eqref{eq:c-hankel} after convolution with $\mathcal K$.
These operations are justified despite the singularity of $\mu$ at zero:
for real $a$, $|1-J_0(ar)|\leq a^2r^2/4$, so
\[
 \int_{\mathbb R}|\mathcal K(w)|\int_0^4h(r)
       |1-J_0((v-w)r)|\,dr\,dw
 \leq\frac14\int_{\mathbb R}|\mathcal K(w)|(v-w)^2\,dw<\infty.
\]
For each $r\leq3$ the inner convolution is exactly zero, since
$\phi(0)=0$ and $\phi$ vanishes on $[-r,r]$.  On $r\geq3$ the density
$h$ is integrable, and a further application of Fubini gives the exact
identity
\begin{equation}\label{eq:c-filtered}
 (\mathcal K*f)(v)=\int_3^4 e^{ivx}\phi(x)q(x)\,dx,
 \qquad
 q(x)=-\frac1\pi\int_x^4\frac{h(r)}{\sqrt{r^2-x^2}}\,dr.
\end{equation}
Thus no unestimated small-$u$ term remains in the filtered integral.

For $y=4-x$, use \eqref{eq:c-hendpoint} and set $r=x+yz$ to get
\[
 q(4-y)=-\frac{y^{5/2}}\pi
 \int_0^1\frac{(1-z)^2z^{-1/2}H(y(1-z))}
                   {\sqrt{8-y(2-z)}}\,dz
 =y^{5/2}B(y),
\]
where $B$ is real analytic at zero, with
\[
 B(0)=-\frac{1}{64\pi^2\sqrt8}
                 B\!\left(\frac12,3\right)
       =-\frac{1}{120\sqrt2\pi^2}.
\]
The Fourier form of Watson's lemma
\cite[Eq.~(2.3.20)]{DLMF}, applied at $x=4$, now yields
\begin{align}\label{eq:c-endpointvalue}
 (\mathcal K*f)(v)
 &=B(0)\Gamma(7/2)e^{-7\pi i/4}
          e^{4iv}v^{-7/2}+\O(v^{-9/2})\notag\\
 &=-\frac{e^{i\pi/4}}{64\sqrt2\pi^{3/2}}
          e^{4iv}v^{-7/2}+\O(v^{-9/2}),\qquad v\to+\infty.
\end{align}
Here the hypotheses of that lemma include differentiable endpoint
expansions; they hold because $B$ is analytic.  Away from $x=4$ the
cutoff integrand is smooth and compactly supported, so repeated
integration by parts makes its contribution $\O(v^{-M})$ for any $M$.

\emph{Identification with the oscillatory BT coefficient.}
The full asymptotic expansions of the actual solutions in\eqref{eq:fullexp}, truncated through order $n^{-3/2}$ in the
algebraic branches and at the leading oscillatory terms, give
\begin{equation}\label{eq:c-BTtruncated}
 F_n=Q_0(n)+n^{-7/4}
       \bigl(K_3\,e^{8i\sqrt n}+\bar K_3 e^{-8i\sqrt n}\bigr)
       +\O(n^{-9/4}),
\end{equation}
where
\[
 Q_0(n)=K_1\sqrt n-\frac12
       +(a\log n+b)n^{-1/2}+(d\log n+e)n^{-3/2}
\]
for real constants $a,b,d,e$; their values will not be needed.  Values of
these expressions at $n=0$ may be assigned arbitrarily.

Let $N$ have the Poisson distribution of mean $t$.  Taylor expansion
on $|N-t|\leq t/2$, together with the exponential Poisson tail bound,
shows that
\begin{equation}\label{eq:c-symbolpoisson}
 \mathbb E Q_0(N)=\widetilde Q(2\sqrt t)
                    +\O(t^{-5/2}\log t),
\end{equation}
where $\widetilde Q(v)$ is a finite linear combination of
$v,1,v^{-1}\log v,v^{-1},v^{-3}\log v,v^{-3}$.
For completeness, for the term $n^{1/2}$ one Taylor expands to degree
$5$, using $\mathbb E(N-t)=0$, $\mathbb E(N-t)^2=t$,
$\mathbb E(N-t)^3=t$, $\mathbb E(N-t)^4=3t^2+t$,
$\mathbb E(N-t)^5=10t^2+t$ and
$\mathbb E|N-t|^6=\O(t^3)$; the remainder is $\O(t^{-5/2})$.
Expansion to degrees $3$ and $1$ suffices for the terms of order
$n^{-1/2}\log n$ and $n^{-3/2}\log n$, respectively, with remainder
$\O(t^{-5/2}\log t)$.  The terms without logarithms are handled in
the same way.

To treat the oscillation, put $Z=(N-t)/\sqrt t$.  On the same central
event,
\[
 8(\sqrt N-\sqrt t)=4Z+\O(t^{-1/2}Z^2),\qquad
 (N/t)^{-7/4}=1+\O(t^{-1/2}|Z|).
\]
Since $\mathbb E Z^2=1$ and
\[
 \mathbb E e^{4iZ}
 =\exp\!\left(t\bigl(e^{4i/\sqrt t}-1-4i/\sqrt t\bigr)\right)
 =e^{-8}\bigl(1+\O(t^{-1/2})\bigr),
\]
we obtain, with the negligible $N=0$ event omitted,
\begin{equation}\label{eq:c-poissonfactor}
 \mathbb E\bigl[N^{-7/4}e^{8i\sqrt N};N\geq1\bigr]
 =t^{-7/4}e^{8i\sqrt t}
       \bigl(e^{-8}+\O(t^{-1/2})\bigr).
\end{equation}
The conjugate formula also holds, and Poissonization of the remainder
in \eqref{eq:c-BTtruncated} is $\O(t^{-9/4})$.
Combining \eqref{eq:c-symbolpoisson}--\eqref{eq:c-poissonfactor} gives
\begin{equation}\label{eq:c-fasymptotic}
 f(v)=\widetilde Q(v)+2^{7/2}e^{-8}v^{-7/2}
                \bigl(K_3\,e^{4iv}+\bar K_3 e^{-4iv}\bigr)
                 +\O(v^{-9/2}),\qquad v\to+\infty.
\end{equation}

Convolution with $\mathcal K$ eliminates $\widetilde Q$ to every algebraic
order.  Indeed, extend its elementary terms smoothly to the real line
with polynomial growth, Taylor expand at $v$ for $|w|\leq v/2$, and
use all the vanishing moments in \eqref{eq:c-filtermoments}; the
Schwartz tails on $|w|>v/2$ give arbitrarily small algebraic errors.
For $v^{-7/2}e^{\pm4iv}$ we likewise first use a fixed smooth cutoff to
extend these functions to the whole real line with polynomial growth;
the same argument then uses $\phi(4)=1$ and $\phi(-4)=0$.
Finally, convolution preserves the
$\O(v^{-9/2})$ remainder: on $|w|\leq v/2$ use its pointwise bound,
and on the complement use its global polynomial bound and the rapid
decay of $\mathcal K$.  All required
global growth bounds hold, since $F_n\leq\sqrt{2n}$ implies
$0\leq f(v)\leq|v|/\sqrt2$.
It follows that
\begin{equation}\label{eq:c-filterBT}
 (\mathcal K*f)(v)=2^{7/2}e^{-8}K_3\,e^{4iv}v^{-7/2}
                    +\O(v^{-9/2}).
\end{equation}
No differentiability of a BT remainder has been used here.
Comparing \eqref{eq:c-endpointvalue} and \eqref{eq:c-filterBT} proves
\[
 2^{7/2}e^{-8}K_3=-\frac{e^{i\pi/4}}{64\sqrt2\pi^{3/2}},
 \qquad K_3=-\frac{e^{8+i\pi/4}}{1024\pi^{3/2}}.
\]
The asserted additive cosine expansion follows from
$Y(n)=e^{8i\sqrt n}n^{-7/4}(1+\O(n^{-1/2}))$ and conjugation.
\end{proof}

\subsection{The connection constants for every hook length.}\label{ssec:consth}
We now evaluate the constants $K_2^{(h)}$ and $K_3^{(h)}$ of Theorem~\ref{mth:hookh}(a); for $h=1$ they are the values of Theorem~\ref{mth:full}(b). The idea is to transfer the expansion of $F_n$ to $F_n^{(h)}$ through an exact finite identity. For an integer $j\geq1$ write $(a)_j=a(a+1)\cdots(a+j-1)$, and let
$$
T_j[f](n)=\frac1{(n+1)_j}\sum_{k=1}^{n}(k+1)_{j-1}\,f(k)
$$
for a sequence $f$.

\begin{proof}[Proof of Theorem~\ref{mth:hookh}\textup{(b)} for $K_2^{(h)}$ and $K_3^{(h)}$]
\emph{An exact identity.} Insert the binomial thinning \eqref{eq:thinning} into \eqref{eq:FhG}. Since
$$
\binom nk\int_0^1s^{k+j-1}(1-s)^{n-k}\,ds=\binom nk\frac{(k+j-1)!\,(n-k)!}{(n+j)!}=\frac{(k+1)_{j-1}}{(n+1)_j},
$$
and $F_0=0$, we get, for every $n\geq1$,
\begin{equation}\label{eq:FhT}
F_n^{(h)}=(-1)^{h-1}h^2\Bigl(F_n+\sum_{j=2}^{h}\beta_j\,T_j[F](n)\Bigr).
\end{equation}

\emph{Sums of powers.} Let $j\geq2$. The Euler--Maclaurin formula, applied to $(k+1)_{j-1}f(k)$ with $(k+1)_{j-1}=k^{j-1}+\frac{j(j-1)}2k^{j-2}+\cdots$, and the expansion $(n+1)_j=n^j\bigl(1+\frac{j(j+1)}{2n}+\O(n^{-2})\bigr)$ give
\begin{equation}\label{eq:Tjpowers}
\begin{aligned}
T_j[1](n)&=\frac1j-\frac{(j-1)!}{(n+1)_j},\\
T_j[\sqrt k\,](n)&=\frac{\sqrt n}{j+\frac12}-\frac{1}{2(4j^2-1)\sqrt n}+\O(n^{-3/2}),\\
T_j[k^{-1/2}](n)&=\frac{1}{(j-\frac12)\sqrt n}+\O(n^{-3/2}),\\
T_j[k^{-1/2}\log k](n)&=\frac{1}{(j-\frac12)\sqrt n}\Bigl(\log n-\frac{1}{j-\frac12}\Bigr)+\O\Bigl(\frac{\log n}{n^{3/2}}\Bigr).
\end{aligned}
\end{equation}
Here the first identity is exact, since $\sum_{k=0}^n(k+1)_{j-1}=(n+1)_j/j$ and the omitted term $k=0$ is $(1)_{j-1}=(j-1)!$; thus $T_j[1](n)=\frac1j+\O(n^{-j})$, and the correction $\O(n^{-j})=\O(n^{-2})$ is harmless below. The constant produced by the Euler--Maclaurin formula contributes $\O(n^{-j})=\O(n^{-2})$. For the second one, the coefficient of $n^{-1/2}$ is $\frac12+\frac{j(j-1)}{2j-1}-\frac{j(j+1)}{2j+1}=-\frac1{2(4j^2-1)}$. Moreover, if $f(k)=\O(k^{-3/2}\log k)$, then $T_j[f](n)=\O(n^{-3/2}\log n)$.

\emph{The constant $K_2^{(h)}$.} By \eqref{eq:shape}, $F_k=K_1\sqrt k+\frac{3K_1}{2^8}\frac{\log k}{\sqrt k}+\frac{K_2}{\sqrt k}-\frac12+\O(k^{-3/2}\log k)$. Insert this in \eqref{eq:FhT} and use \eqref{eq:Tjpowers}. With $R(m)=\rho_h(m)/h^2=1+\sum_j\beta_j/(m+j)$, we have
$$
1+\sum_j\frac{\beta_j}{j+\frac12}=R(\tfrac12),\quad
1+\sum_j\frac{\beta_j}{j-\frac12}=R(-\tfrac12),
$$
$$
\sum_j\frac{\beta_j}{(j-\frac12)^2}=-R'(-\tfrac12),\quad
\sum_j\frac{\beta_j}{2(4j^2-1)}=\frac{R(-\frac12)-R(\frac12)}{8}.
$$
Hence the coefficient of $n^{-1/2}$ in $F_n^{(h)}$ is
$$
(-1)^{h-1}h^2\Bigl(K_2R(-\tfrac12)+\frac{3K_1}{2^8}R'(-\tfrac12)-\frac{K_1}8\bigl(R(-\tfrac12)-R(\tfrac12)\bigr)\Bigr).
$$
From the product formula for $\rho_h$ we get $R(-\frac12)=(-1)^{h-1}$, $R(\frac12)=(-1)^{h-1}\frac{3}{4h^2-1}$ (as in the proof of \eqref{eq:asyh}), and
$$
\frac{R'(-\frac12)}{R(-\frac12)}=\sum_{i=1}^{h-1}\frac{1}{-\frac12-i}-\sum_{i=2}^{h}\frac{1}{i-\frac12}=-4\sum_{i=1}^{h-1}\frac1{2i+1}.
$$
The same computation recovers the coefficients $K_1^{(h)}$, $-\frac h2$ and $\frac{h^2}{8\pi^2}$ of $\sqrt n$, $1$ and $\frac{\log n}{\sqrt n}$. Since $S_1^{(h)}$ has no term in $n^{-1/2}$, part~(a) shows that $K_2^{(h)}$ is the coefficient of $n^{-1/2}$ in $F_n^{(h)}$. Therefore
$$
K_2^{(h)}=h^2\Bigl(K_2-\frac{3K_1}{64}\sum_{i=1}^{h-1}\frac1{2i+1}-\frac{K_1(h^2-1)}{2(4h^2-1)}\Bigr).
$$
With $K_1=\frac{32}{3\pi^2}$, the value of $K_2$ from Theorem~\ref{mth:full}(b), and $\sum_{i=1}^{h-1}\frac1{2i+1}=\sum_{i=1}^h\frac1{2i-1}-1$, this is the value stated in part~(b).

\emph{The constant $K_3^{(h)}$.} By \eqref{eq:fullexp} and Lemma~\ref{lem:BT}, for $k$ large enough,
$$
F_k=\Phi(k)+K_3Y(k)+\bar K_3\bar Y(k)+\O(k^{-5/2}\log k),
$$
where $\Phi(k)$ is the sum of the terms of order at least $k^{-3/2}$ in the expansion \eqref{eq:fullexpM}; the finitely many small $k$ contribute $\O(n^{-j})$ to $T_j[F](n)$. Let $j\geq2$ and $f(x)=x^{j-11/4}e^{8i\sqrt x}$. The substitution $x=t^2$ and an integration by parts give
$$
\int_1^nf(x)\,dx=2\int_1^{\sqrt n}t^{2j-9/2}e^{8it}\,dt=\Bigl[\frac{t^{2j-9/2}e^{8it}}{4i}\Bigr]_1^{\sqrt n}-\frac{2j-9/2}{4i}\int_1^{\sqrt n}t^{2j-11/2}e^{8it}\,dt .
$$
For $j\geq3$ the last integral is $\O(n^{j-9/4})$; for $j=2$ it converges as $n\to\infty$, with a tail $\O(n^{-1/4})$. In both cases $\int_1^nf(x)\,dx=C_0+\O(n^{j-9/4})$ for a constant $C_0$. Since $f'(x)=\O(x^{j-13/4})$, the difference between $\sum_{k=1}^nf(k)$ and $\int_1^nf(x)\,dx$ is $\frac{f(1)+f(n)}2$ plus an integral bounded by $\int_1^n|f'(x)|\,dx$, which is finite for $j=2$ and $\O(n^{j-9/4})$ for $j\ge3$ (Euler--Maclaurin); hence $\sum_{k\le n}f(k)=C+\O(n^{j-9/4})$ for a constant $C$. Moreover, $(k+1)_{j-1}Y(k)$ differs from its leading oscillatory term $f(k)$ by $\O(k^{j-13/4})$; summing this error absolutely gives again a constant plus $\O(n^{j-9/4})$. Dividing by $(n+1)_j$, we get $T_j[Y](n)=\O(n^{-j})+\O(n^{-9/4})=\O(n^{-2})$, and similarly $T_j[\bar Y](n)=\O(n^{-2})$. Similarly, $T_j$ applied to the remainder $\O(k^{-5/2}\log k)$ is $\O(n^{-2})$. Finally, by the Euler--Maclaurin formula, $T_j[\Phi](n)$ is, up to an error $\O(n^{-2})$, a finite linear combination of $1$, $n^{1/2-m}$ and $n^{1/2-m}\log n$. By \eqref{eq:FhT}, we therefore get
$$
F_n^{(h)}=\Psi(n)+(-1)^{h-1}h^2\bigl(K_3Y(n)+\bar K_3\bar Y(n)\bigr)+\O(n^{-2}),
$$
where $\Psi(n)$ is a finite linear combination of $1$, $n^{1/2-m}$ and $n^{1/2-m}\log n$ with $0\le m\leq2$. Compare this with part~(a). The difference between $\Psi(n)$ and the non-oscillatory part $K_1^{(h)}S_1^{(h)}+K_2^{(h)}S_2^{(h)}-\frac h2$ is $\O(n^{-7/4})$ and has an expansion in the same powers and logarithms. Hence its coefficients vanish up to order $n^{-3/2}$, and it is $\O(n^{-5/2}\log n)$. Since $Y(n)$ and $Y^{(h)}(n)$ are both equal to $e^{8i\sqrt n}n^{-7/4}(1+\O(n^{-1/2}))$, we get
$$
2\Re\Bigl(\bigl(K_3^{(h)}-(-1)^{h-1}h^2K_3\bigr)e^{8i\sqrt n}\Bigr)n^{-7/4}=\O(n^{-2}).
$$
As in the proof of part~(a), the phases $8\sqrt n$ are dense modulo $2\pi$, so $K_3^{(h)}=(-1)^{h-1}h^2K_3$.
\end{proof}

Thus the oscillatory part of $F_n^{(h)}$ is $(-1)^{h-1}h^2$ times that of $F_n$, up to $\O(n^{-9/4})$. This is consistent with the localized proof in Section~\ref{ssec:cconst}: the oscillation comes from the endpoint $u=16$, that is, from the behavior of the coefficients for large $m$, where $(-1)^{h-1}\rho_h(m)\sim(-1)^{h-1}h^2$.

\subsection{The expected sum of the hook lengths.}\label{ssec:hooksum}
Let
$$
u_n:=\mathbb{E}\Bigl[\sum_{\square\in\la_{(n)}^{(2)}}h_\square\Bigr].
$$
Schachinger \cite[Section~2.1]{Schachinger2023} showed that
\begin{equation}\label{eq:hsun}
u_n=\sum_{m=1}^n\binom nm b_m,\qquad b_m=\frac{(-1)^m}{(2m^2-3m)\,m!}\binom{2m-2}{m-1}^{2},
\end{equation}
so that $u_0,\dots,u_4=0,1,3,\frac{16}3,\frac{49}6$. Recall that $F_n=\sum_m\binom nma_m$ with $a_m=(-1)^{m-1}C_{m-1}C_m/m!$ (Theorem~\ref{th:exactF2_sec1}).

\begin{prop}\label{prop:hsidentity}
Let $\varrho(m)=-8-\frac6m+\frac{30}{m+1}=-\frac{2(2m-1)(2m-3)}{m(m+1)}$. Then $a_m=\varrho(m)b_m$ for $m\geq1$, and for every $n\geq0$,
\begin{equation}\label{eq:hsidentity}
F_n=-8u_n-6\sum_{k=1}^n\frac{u_k}{k}+\frac{30}{n+1}\sum_{k=0}^nu_k .
\end{equation}
Conversely, with $P_c(n)=\Gamma(n+1)/\Gamma(n+1-c)$,
\begin{equation}\label{eq:hsinverse}
u_n=-\frac{F_n}8+\frac3{32}\,P_{1/2}(n)\sum_{k=1}^n\frac{\Gamma(k-\frac12)}{k!}F_k-\frac{15}{32}\,P_{3/2}(n)\sum_{k=1}^n\frac{\Gamma(k-\frac32)}{k!}F_k .
\end{equation}
\end{prop}

\begin{proof}
From $C_{m-1}=\frac1m\binom{2m-2}{m-1}$ and $C_m=\frac{2(2m-1)}{m+1}C_{m-1}$ we get $a_m/b_m=\varrho(m)$. Since $\sum_{k=m}^n\frac1k\binom km=\frac1m\binom nm$ and $\sum_{k=m}^n\binom km=\frac{n+1}{m+1}\binom nm$, the two sums in \eqref{eq:hsidentity} are $\sum_m\binom nm\frac{b_m}m$ and $(n+1)\sum_m\binom nm\frac{b_m}{m+1}$, which proves \eqref{eq:hsidentity}. Conversely, $\frac1{\varrho(m)}=-\frac18+\frac{3}{32}\frac1{m-1/2}-\frac{15}{32}\frac1{m-3/2}$. For $c\in\{\frac12,\frac32\}$, the sequence $G^{(c)}_n=\sum_m\binom nm\frac{a_m}{m-c}$ satisfies $(n-c)G^{(c)}_n=nG^{(c)}_{n-1}+F_n$, because $n\binom{n-1}{m-1}=m\binom nm$; with $G^{(c)}_0=0$ this gives $G^{(c)}_n=P_c(n)\sum_{k=1}^nF_k\Gamma(k-c)/k!$, and \eqref{eq:hsinverse} follows.
\end{proof}

The sequence $u_n$ satisfies the homogeneous recurrence
\begin{equation}\label{eq:hsrec}
\begin{aligned}
&(n+1)(n+2)(n+3)u_n-(4n-7)(n+2)(n+3)u_{n+1}+(2n-1)(3n+2)(n+3)u_{n+2}\\
&\qquad-(4n^3+23n^2+63n+78)u_{n+3}+(n+3)(n+4)^2u_{n+4}=0
\end{aligned}
\end{equation}
for $n\geq0$; it is stated in \cite[Section~2.1]{Schachinger2023}. A certificate, as in the proof of Proposition~\ref{prop:rec}, is $G_{n,m}=-\frac{m^3(m-1)}{n+4}\binom{n+4}mb_m$. Its formal solutions are found as in Proposition~\ref{prop:formal}. The indicial polynomial is $4(2\theta-1)(2\theta-3)$, and there are four formal solutions $U_1,U_2,W,\bar W$: $U_2$ is a series in $n^{1/2-m}$; $U_1=A+\frac{27}{256}\,U_2\log n$, where $A$ is a series in $n^{3/2-m}$ without a term in $n^{1/2}$; and $W$ is $e^{8i\sqrt n}n^{-7/4}$ times a series in $n^{-1/2}$. The first terms are
\begin{align*}
A(n)&=n^{3/2}-\frac{1559}{2^{18}}n^{-1/2}-\frac{131511}{2^{26}}n^{-3/2}+\frac{196021423}{2^{38}}n^{-5/2}+\cdots,\\
U_2(n)&=n^{1/2}+\frac{65}{2^9}n^{-1/2}+\frac{1855}{2^{18}}n^{-3/2}-\frac{1206355}{2^{28}}n^{-5/2}+\cdots,\\
W(n)&=e^{8i\sqrt n}n^{-7/4}\Bigl(1+\frac{221\,i}{192\,n^{1/2}}+\frac{13871}{73728\,n}+\cdots\Bigr).
\end{align*}

\begin{cor}\label{cor:hooksum}
For $n$ large enough,
$$
u_n=\kappa_1U_1(n)+\kappa_2U_2(n)+\kappa_3W(n)+\bar\kappa_3\bar W(n),
$$
where $U_1,U_2,W,\bar W$ are the solutions of \eqref{eq:hsrec} with the above asymptotic expansions, and
$$
\kappa_1=\frac{256}{27\pi^2},\qquad\kappa_2=\frac{2\gamma+12\log2-\frac{77}9}{\pi^2},\qquad\kappa_3=-\frac{K_3}8=\frac{e^{8+i\pi/4}}{8192\,\pi^{3/2}}.
$$
Equivalently, with $\delta_n=\log n+2\gamma+12\log2$, the non-oscillatory part of $u_n$ is $\frac{256}{27\pi^2}A(n)+\frac{\delta_n-77/9}{\pi^2}\,U_2(n)$, to all orders. In particular,
\begin{align*}
u_n={}&\frac{256}{27\pi^2}n^{3/2}+\frac{9\delta_n-77}{9\pi^2}n^{1/2}
+\frac{3510\delta_n-31589}{27648\pi^2}n^{-1/2}+\frac{5565\delta_n-62224}{786432\pi^2}n^{-3/2}\\
&+\frac{e^8}{2^{12}\pi^{3/2}}\,n^{-7/4}\Bigl(\cos\theta_n-\frac{221}{192\sqrt n}\sin\theta_n\Bigr)
+\frac{1310693443-130286340\,\delta_n}{27\cdot2^{30}\,\pi^2}\,n^{-5/2}+\O\bigl(n^{-11/4}\bigr),
\end{align*}
where $\theta_n=8\sqrt n+\frac\pi4$.
\end{cor}

The expansion up to $\O(n^{-9/4})$ is \cite[Theorem~2.1]{Schachinger2023}, there stated for $\mathbb{E}[\sum_\square h_\square]-n$. Corollary~\ref{cor:hooksum} gives an independent proof of it, together with all the higher-order terms.

\begin{proof}
The four formal solutions are formally independent, so Lemma~\ref{lem:BT} gives the decomposition, with $\kappa_1,\kappa_2$ real; in particular $u_n=\kappa_1n^{3/2}+\frac{27\kappa_1}{256}n^{1/2}\log n+\kappa_2n^{1/2}+\O(n^{-1/2}\log n)$.

\emph{The constant $\kappa_3$.} Apply \eqref{eq:hsidentity} to the decomposition. The partial sums of $U_1(k)/k$ and $U_1(k)$, and those of $U_2$, have expansions in powers of $n$ and $\log n$, up to a constant and an error $\O(n^{-5/2}\log n)$ (Euler--Maclaurin). For the oscillatory solutions, $\sum_{k>n}e^{8i\sqrt k}k^{-a}=\O(n^{1/2-a})$ for $a>1$, by comparison with the integral and one integration by parts, as in the proof of $K_3^{(h)}$ in Section~\ref{ssec:consth}. Hence $\sum_{k\le n}W(k)/k=C+\O(n^{-9/4})$ and $\frac1{n+1}\sum_{k\le n}W(k)=\frac{C'}{n+1}+\O(n^{-9/4})$. So the oscillatory part of the right side of \eqref{eq:hsidentity} is $-8\bigl(\kappa_3e^{8i\sqrt n}+\bar\kappa_3e^{-8i\sqrt n}\bigr)n^{-7/4}+\O(n^{-9/4})$. Comparing with Theorem~\ref{mth:full}, and using that the phases $8\sqrt n$ are dense modulo $2\pi$, gives $-8\kappa_3=K_3$.

\emph{The constants $\kappa_1$ and $\kappa_2$.} By \eqref{eq:k2moments}, $C_{m-1}C_m=M(m)$, so $b_m=(-1)^{m-1}N(m)/m!$ with $N(r)=M(r)/\varrho(r)$. For $c\in\{\frac12,\frac32\}$, the nonnegative function $t^{-c-1}\int_t^{16}u^c\mu(u)\,du$ has $r$-th moment $M(r)/(r-c)$ for $r>\frac32$. By the partial fractions of $1/\varrho$, $N(r)$ is therefore the $r$-th moment of a combination of three nonnegative densities, for $r>\frac32$. Let $\beta(q)=\sum_{m\ge1}b_mq^m$ (note $b_1=1$). The argument leading to \eqref{eq:k2mellin}, with $e^{-x}-1+x$ in place of $1-e^{-x}$, shows that $\beta(q)-q$ has the Mellin transform
$$
-\Gamma(s)N(-s)=\frac{\Phi(s)}{\varrho(-s)}=\frac{16^{-s}}{32\pi}\,\frac{\Gamma(s)\Gamma(-s-\frac32)\Gamma(-s-\frac12)}{\Gamma(-s)\Gamma(1-s)}
$$
in the strip $-2<\Re s<-\frac32$. It decays like $\Phi$ on vertical lines. Shift the line of integration to $\Re s=0$, where the integrand is regular. The poles crossed are $s=-\frac32$ (simple), $s=-1$ (simple, with residue $q$, which cancels the term $q$) and $s=-\frac12$ (double). This gives
$$
\beta(q)=\alpha_1q^{3/2}+q^{1/2}(\alpha_2\log q+\alpha_3)+\O(1),\quad
\alpha_1=\frac{64}{9\pi^{3/2}},\ \alpha_2=\frac1{2\pi^{3/2}},\ \alpha_3=\frac{\gamma+10\log2-3}{2\pi^{3/2}} .
$$
Since $\sum_{n\ge0}u_nz^n=\beta(z/(1-z))/(1-z)$, the Abelian comparison of Section~\ref{ssec:k2evaluation}, now with $\sum_nn^{3/2}e^{-tn}=\Gamma(\frac52)t^{-5/2}+\O(t^{-1/2})$, gives the three relations
$$
\kappa_1\Gamma(\tfrac52)=\alpha_1,\qquad
\tfrac{27}{256}\kappa_1\Gamma(\tfrac32)=\alpha_2,\qquad
\kappa_2\Gamma(\tfrac32)=\alpha_3-\frac{\alpha_1}4-\alpha_2\psi(\tfrac32).
$$
With $\psi(\frac32)=2-\gamma-2\log2$, they give the stated $\kappa_1$ and $\kappa_2$.
\end{proof}

\smallskip\noindent\emph{Remark.} By \eqref{eq:hsinverse}, the complete expansion of $u_n$ also follows directly from that of $F_n$, except for the coefficients of the two homogeneous solutions $P_{3/2}(n)\sim n^{3/2}$ and $P_{1/2}(n)\sim n^{1/2}$. These coefficients depend on global quantities, and not on the asymptotics of $F_n$: for $c=\frac32$, on the convergent sum $\sum_{k\ge1}F_k\Gamma(k-\frac32)/k!$; for $c=\frac12$ the sum diverges, since $F_k\Gamma(k-\frac12)/k!=K_1/k+\O(k^{-3/2})$, and one uses instead its finite part $\lim_{n\to\infty}\bigl(\sum_{k\le n}F_k\Gamma(k-\frac12)/k!-K_1\sum_{k\le n}\frac1k\bigr)$, the harmonic sum $\sum_{k\le n}\frac1k=\log n+\gamma+\O(n^{-1})$ contributing the term in $P_{1/2}(n)\log n$ and an extra $K_1\gamma$ to the coefficient of $P_{1/2}(n)$. These coefficients correspond to the zeros $\varrho(\frac32)=\varrho(\frac12)=0$. This is the mirror image of the fact that the constant $-\frac12$ in $F_n$ is not determined by the asymptotics of $u_n$.

\section{Discussion and future work}\label{sec:discussion}

\subsection{The uniform and Gelfand measures.}\label{ssec:beta01}
For $\beta=0$ and $\beta=1$, the complete asymptotic expansion of $\mathbb{E}[\co(\la_{(n)}^{(\beta)})]$ is already known.

For the uniform measure ($\beta=0$), the generating function of $\sum_{|\la|=n}\co(\la)$ is $\frac{x}{1-x}\prod_{k\geq1}\frac1{1-x^k}$, and Grabner, Knopfmacher and Wagner \cite[Proposition~4.1]{Grabner2014KW} derived from it
\begin{equation}\label{eq:beta0exp}
\mathbb{E}[\co(\la_{(n)}^{(0)})]=\frac{\sqrt6}{\pi}\sqrt n+\frac{6-\pi^2}{2\pi^2}+\frac{2\pi^4-3\pi^2+216}{24\pi^3\sqrt{6n}}+\frac{54-\pi^4}{12\pi^4n}+\O\bigl(n^{-3/2}\bigr);
\end{equation}
their method gives arbitrarily many terms \cite[Remark~4]{Grabner2014KW}.

For the Gelfand measure ($\beta=1$), the problem reduces to fixed points of involutions. Let $\Invo_n$ be the set of involutions of $[n]$, let $\fix(\sigma)$ be the number of fixed points of $\sigma$, and let $a_n=\#\Invo_n$.

\begin{thm}\label{th:fix=cor_sec1}
For every positive integer $n$,
$$
\sum_{|\la|=n}f_\la\,\co(\la)=\sum_{\sigma\in\Invo_n}\fix(\sigma).
$$
\end{thm}

\begin{proof}
The left-hand side counts the pairs $(T,c)$, where $T$ is a standard Young tableau with $n$ boxes and $c$ is a corner of $T$. The right-hand side counts the pairs $(\sigma,i)$, where $\sigma\in\Invo_n$ and $i$ is a fixed point of $\sigma$. Given $(T,c)$, apply reverse row insertion (reverse bumping) to $T$ starting from the corner $c$ \cite[Section~7.11]{Stanley1999EC2}. It produces a letter $i\in[n]$ and a standard tableau $S$ of shape $\shape(T)\setminus\{c\}$ on the ordered alphabet $[n]\setminus\{i\}$. By the Robinson--Schensted correspondence, $S$ is the insertion tableau of a unique involution $\sigma'$ of $[n]\setminus\{i\}$ (namely, $(S,S)$ is the image of $\sigma'$). Let $\sigma$ be the involution of $[n]$ obtained by adding the fixed point $i$ to $\sigma'$, and map $(T,c)$ to $(\sigma,i)$. Conversely, given $(\sigma,i)$, remove the fixed point $i$ from $\sigma$, take the insertion tableau $S$ of the resulting involution of $[n]\setminus\{i\}$, row-insert $i$ into $S$ to get $T$, and let $c$ be the new box. Since row insertion and reverse row insertion are inverse to each other, these two maps are inverse bijections. Note that the letter $i$ is in general neither $n$ nor the entry of $T$ in the corner $c$.
\end{proof}

By the Robinson--Schensted correspondence, $\sum_{|\mu|=n}f_\mu=a_n$. Hence, by Theorem~\ref{th:fix=cor_sec1}, $\mathbb{E}[\co(\la_{(n)}^{(1)})]$ is the mean number of fixed points of a uniform random involution of $[n]$; since each $k\in[n]$ is a fixed point of exactly $a_{n-1}$ involutions, it is equal to $n\,a_{n-1}/a_n$. Lugo \cite[p.~188]{Lugo2010} obtained the complete asymptotic expansion of this mean, by inserting the asymptotic series of $a_n$ into the recurrence $a_n=a_{n-1}+(n-1)a_{n-2}$:
\begin{equation}\label{eq:beta1exp}
\mathbb{E}[\co(\la_{(n)}^{(1)})]=\sqrt n-\frac12+\frac{3}{8\sqrt n}-\frac1{8n}-\frac{1}{128\,n^{3/2}}+\frac{3}{32\,n^2}-\frac{85}{1024\,n^{5/2}}+\O\bigl(n^{-3}\bigr).
\end{equation}
The constant term $-\frac12$ is the same as for $\beta=2$ (Theorem~\ref{mth:full}), whereas it is $\frac3{\pi^2}-\frac12$ for $\beta=0$.

\subsection{The leading constant $c(\beta)$.}
By \eqref{eq:beta0exp}, \eqref{eq:beta1exp} and Theorem~\ref{mth:full}, the leading constants are
$$
c(0)=\frac{\sqrt6}{\pi}=0.7797\ldots,\qquad c(1)=1,\qquad c(2)=\frac{32}{3\pi^2}=1.0808\ldots,
$$
which are strictly increasing in $\beta$. These results make us conjecture the following asymptotic behavior of general $\mathbb{E}[\co(\la_{(n)}^{(\beta)})]$.

\begin{con}\label{con:cbeta}
For each given $\beta\geq 0$,
when $n$ tends to infinity, we have
	\begin{equation}
		\mathbb{E}[\co(\la_{(n)}^{(\beta)})] \sim c(\beta) \sqrt{n}
	\end{equation}
	where $c(\beta)>0$ is a constant determined by $\beta$. Furthermore, $c(\beta)$ is a non-decreasing function of $\beta$, i.e.,  $c(\beta_1)\leq c(\beta_2)$ whenever $\beta_1 \leq \beta_2$. 
\end{con}

\subsection{An upper bound for $c(\beta)$.}\label{ssec:cupper}
Lemma~\ref{lem: upper_bound} gives the trivial bound $c(\beta)\le\sqrt2=1.4142\ldots$. For $\beta>0$ the limit shape yields a better bound. Encode a partition $\la$ by its Maya set $S(\la)=\{\la_i-i: i\geq 1\}\subset\mathbb Z$ (with $\la_i=0$ for $i>\ell(\la)$), whose elements we call particles and whose complement we call holes. A row of $\la$ ending a part of a new size corresponds to a particle $x$ with $x-1$ a hole, so that
\begin{equation}\label{eq:cornerMaya}
\co(\la)=\#\{x\in\mathbb Z:\ x\in S(\la),\ x-1\notin S(\la)\}.
\end{equation}
Under the Vershik--Kerov--Logan--Shepp limit shape $\Omega$, the particles near position $x=u\sqrt n$ have density
$$
\rho(u)=\frac{1-\Omega'(u)}{2}=\frac1\pi\arccos\frac u2\qquad(-2\le u\le 2),
$$
with $\rho=1$ for $u<-2$ and $\rho=0$ for $u>2$.

\begin{proof}[Proof of Theorem~\ref{mth:beta}]
We prove that, for every $\beta>0$,
\begin{equation}\label{eq:cupper}
\limsup_{n\to\infty}\frac{\mathbb{E}[\co(\la_{(n)}^{(\beta)})]}{\sqrt n}\le\int_{-2}^{2}\min\bigl(\rho(u),1-\rho(u)\bigr)\,du=\frac4\pi=1.2732\ldots
\end{equation}

(i) \emph{The Durfee bound.} Let $D(\la)$ be the side of the Durfee square of $\la$, i.e., the largest $D$ with $\la_D\geq D$, so that $\la_{D+1}\le D$ and $D(\la)^2\le n$. The first $D$ rows of $\la$ contribute at most $D$ distinct part sizes, and the other rows have length at most $\la_{D+1}\le D$, so they contribute at most $D$ distinct part sizes as well. Hence, for every partition $\la$,
\begin{equation}\label{eq:durfee}
\co(\la)\le 2D(\la).
\end{equation}

(ii) \emph{The limit shape for every $\beta>0$.} Let $\omega_\la=\omega(\la)$ be the (rescaled) shape of $\la$ defined in Section~1.1, whose limit for $\beta=2$ is $\Omega$, and let $\|\cdot\|_\infty$ be the sup norm. We show that $\mathbb{P}_n^{(\beta)}\bigl(\|\omega_\la-\Omega\|_\infty>\varepsilon\bigr)\to0$ for every $\varepsilon>0$. Two remarks will be used repeatedly: $\sum_{\la\vdash n}(f_\la)^2=n!$ gives $(f_\la)^2\le n!$ for every $\la$ and $\max_{\la\vdash n}(f_\la)^2\ge n!/p(n)$; hence, for every set $A$ of partitions of $n$,
\begin{equation}\label{eq:Pbeta}
\mathbb{P}_n^{(\beta)}(A)\le p(n)^{1+\beta/2}\,\max_{\la\in A}\Bigl(\frac{(f_\la)^2}{n!}\Bigr)^{\beta/2},
\end{equation}
where $p(n)\le e^{\pi\sqrt{2n/3}}$.

\emph{(a) Long rows and columns are unlikely.} Let $\la_1=k$. A standard Young tableau of shape $\la$ is determined by the set of entries of its first row and by the standardization of the remaining entries, which is a standard Young tableau of shape $\bar\la=(\la_2,\la_3,\dots)$. Hence $f_\la\le\binom nk f_{\bar\la}$, and since $(f_{\bar\la})^2\le(n-k)!$,
$$
\frac{(f_\la)^2}{n!}\le\binom nk^2\frac{(n-k)!}{n!}=\frac{1}{k!}\binom nk\le\Bigl(\frac{e^2n}{k^2}\Bigr)^k,
$$
using $\binom nk\le(en/k)^k$ and $k!\ge(k/e)^k$. By \eqref{eq:Pbeta}, for $K>e$,
$$
\mathbb{P}_n^{(\beta)}\bigl(\la_1\ge K\sqrt n\bigr)\le p(n)^{1+\beta/2}\Bigl(\frac{e^2}{K^2}\Bigr)^{\beta K\sqrt n/2}
\le \exp\Bigl(\sqrt n\Bigl(\pi\sqrt{2/3}\,\bigl(1+\tfrac\beta2\bigr)-\beta K\log\frac Ke\Bigr)\Bigr),
$$
which tends to $0$ if $K=K_\beta$ is large enough. Since $f_{\la'}=f_\la$, the same holds for the length $\la'_1$ of the first column. On the event $\la_1,\la'_1\le K\sqrt n$, the function $g=\omega_\la-\Omega$ is supported in $[-K,K]$ (we take $K\geq 2$) and is $2$-Lipschitz.

\emph{(b) The Vershik--Kerov variational formula.} For a function $g$ on $\mathbb R$ put
$$
\|g\|_{1/2}^2=\int\!\!\int\Bigl(\frac{g(s)-g(t)}{s-t}\Bigr)^2ds\,dt ,
$$
and for $\la\vdash n$ let $F_\la(t)=\sqrt n\,\bigl(\omega_\la(t/\sqrt n)-\Omega(t/\sqrt n)\bigr)$, so that $\|F_\la\|_{1/2}^2=n\,\|\omega_\la-\Omega\|_{1/2}^2$. Vershik and Kerov \cite[Lemma~1 and formulas (5), (8), (9)]{Vershik1985Kerov} proved an exact formula for the Plancherel probability $(f_\la)^2/n!$ of a single diagram; in the form stated in \cite[formula (2)]{Bufetov2012}, it reads
$$
-\log\frac{(f_\la)^2}{n!}=\frac18\,\|F_\la\|_{1/2}^2+\sum_{k\geq1}h_k(\la)\,\gamma_k+\int_{|t|\ge2\sqrt n}F_\la(t)\,\mathrm{arccosh}\frac{|t|}{2\sqrt n}\,dt-\sqrt n\,\epsilon_n ,
$$
where $h_k(\la)$ is the number of boxes of $\la$ with hook length $k$, $\gamma_k=\sum_{l\ge1}\bigl(l(l+1)(2l+1)k^{2l}\bigr)^{-1}>0$, and $\epsilon_n$ depends only on $n$ (not on $\la$) and tends to $0$. The second term is nonnegative, and so is the third, because the profile of $\la$ lies above $|t|$, so that $F_\la(t)\geq0$ for $|t|\geq2\sqrt n$. Hence, for every $\la\vdash n$,
\begin{equation}\label{eq:VKbound}
\frac{(f_\la)^2}{n!}\le\exp\Bigl(-\frac n8\,\|\omega_\la-\Omega\|_{1/2}^2+\sqrt n\,\epsilon_n\Bigr).
\end{equation}
(This is the estimate behind the limit shape theorem of Logan--Shepp and Vershik--Kerov \cite{Logan1977Shepp,Vershik1977Kerov}; see also \cite[Chapter~1]{Romik2015}.)

\emph{(c) From $\|\cdot\|_{1/2}$ to $\|\cdot\|_\infty$.} Let $g$ be $2$-Lipschitz and supported in $[-K,K]$, and let $M=\|g\|_\infty=|g(s_0)|$. Then $|g(s)|\geq M/2$ for $|s-s_0|\le M/4$, so $\int g^2\geq M^3/8$. On the other hand, taking $s\in[-K,K]$ and $t\in[K+1,K+2]$, where $g(t)=0$ and $(s-t)^2\le(2K+2)^2$, we get $\|g\|_{1/2}^2\geq(2K+2)^{-2}\int g^2$. Hence $\|g\|_{1/2}^2\geq M^3/\bigl(8(2K+2)^2\bigr)$.

Combining (a), (b), (c) and \eqref{eq:Pbeta}: on the event $\{\la_1,\la_1'\le K\sqrt n,\ \|\omega_\la-\Omega\|_\infty>\varepsilon\}$, the bound \eqref{eq:VKbound} and (c) give
$$
\frac{(f_\la)^2}{n!}\le\exp\Bigl(-\frac{\varepsilon^3}{64(2K+2)^2}\,n+\sqrt n\,\epsilon_n\Bigr),
$$
so that the probability of this event is at most $p(n)^{1+\beta/2}\exp\bigl(-\beta\varepsilon^3n/(128(2K+2)^2)+\beta\sqrt n\,\epsilon_n/2\bigr)\to0$. Thus, for every $\beta>0$, $\omega_\la\to\Omega$ uniformly, in probability.

(iii) \emph{Conclusion.} The corner $(D,D)$ of the Durfee square lies on the boundary of the Young diagram and on its diagonal. After the rotation by $45^\circ$ and the scaling by $\sqrt2/\sqrt n$ of Section~1.1, it becomes the point of abscissa $0$ and height $2D(\la)/\sqrt n$, so that $\omega_\la(0)=2D(\la)/\sqrt n$. By (ii), $2D(\la)/\sqrt n\to\Omega(0)=\frac4\pi$ in probability under $\mathbb{P}_n^{(\beta)}$. Since $0\le 2D(\la)/\sqrt n\le2$, the convergence also holds in expectation, and \eqref{eq:durfee} gives
$$
\limsup_{n\to\infty}\frac{\mathbb{E}[\co(\la_{(n)}^{(\beta)})]}{\sqrt n}\le\lim_{n\to\infty}\frac{\mathbb{E}[2D(\la_{(n)}^{(\beta)})]}{\sqrt n}=\frac4\pi.
$$
Finally, since $\Omega'=1-2\rho$ and $\rho(-u)=1-\rho(u)$, we have $\Omega(0)=2-\int_0^2(1-2\rho)\,du=2\int_0^2\rho\,du=\int_{-2}^0(1-\rho)\,du+\int_0^2\rho\,du=\int_{-2}^2\min(\rho,1-\rho)\,du$; explicitly, with $u=2\cos\pi r$,
$$
\int_{-2}^2\min(\rho,1-\rho)\,du=2\pi\int_0^1\min(r,1-r)\sin(\pi r)\,dr=4\pi\int_0^{1/2}r\sin(\pi r)\,dr=\frac4\pi.
$$
\end{proof}

The bound is consistent with $c(1)=1$ and $c(2)=32/(3\pi^2)=1.0808\ldots$.

\medskip
\noindent{\it Remark (a local picture).} Formula \eqref{eq:cornerMaya} suggests that $c(\beta)$ should be the integral over $u$ of the density $\mathrm{p}_\beta(\rho)$ of the local pattern ``particle at $x$, hole at $x-1$'', for a local configuration of density $\rho=\rho(u)$. For $\beta=2$, the local process is the discrete sine process \cite{Borodin2000OO}, for which $\mathrm{p}_2(\rho)=\rho-\bigl(\rho^2-\frac{\sin^2\pi\rho}{\pi^2}\bigr)=\rho(1-\rho)+\frac{\sin^2\pi\rho}{\pi^2}$. With $u=2\cos\pi r$ one finds
$$
\int_{-2}^2\Bigl(\rho(1-\rho)+\frac{\sin^2\pi\rho}{\pi^2}\Bigr)du
=2\pi\int_0^1\Bigl(r(1-r)+\frac{\sin^2\pi r}{\pi^2}\Bigr)\sin(\pi r)\,dr
=\frac8{\pi^2}+\frac{8}{3\pi^2}=\frac{32}{3\pi^2},
$$
in agreement with \eqref{eq:asyF2}. For $\beta=2$ this local picture is not a heuristic: it is exactly the proof of Bufetov \cite[Lemma~2.1 and formula~(31)]{Bufetov2012} for $k=1$.

For general $\beta$, pad every partition of $n$ by zero parts to a fixed length $N\geq n$ and put $l_i=\la_i+N-i$; the Frobenius formula $f_\la=n!\prod_{i<j}(l_i-l_j)/\prod_il_i!$ shows that $\mathbb{P}_n^{(\beta)}$ is a discrete logarithmic gas with Dyson index $\beta$. Heuristically, as $\beta\to\infty$ the gas crystallizes: locally, every particle is isolated when $\rho<\frac12$ and every hole is isolated when $\rho>\frac12$, so that $\mathrm{p}_\infty(\rho)=\min(\rho,1-\rho)$ and the bound \eqref{eq:cupper} becomes an equality. As $\beta\to0^+$, the particles should become locally independent with density $\rho$, so that $\mathrm{p}_{0^+}(\rho)=\rho(1-\rho)$, whose integral is $2\pi\int_0^1r(1-r)\sin(\pi r)\,dr=\frac8{\pi^2}$. The macroscopic limit shape alone does not imply these local limits, and the limits $n\to\infty$ and $\beta\to0^+$ need not commute: for fixed $n$, $\mathbb{P}_n^{(\beta)}$ tends to the uniform measure as $\beta\to0^+$, so that $\lim_{n\to\infty}\lim_{\beta\to0^+}\mathbb{E}[\co(\la_{(n)}^{(\beta)})]/\sqrt n=\frac{\sqrt6}{\pi}$ by \eqref{eq:beta0exp}. This suggests the following refinement of Conjecture~\ref{con:cbeta}.

\begin{con}\label{con:cinfty}
We have
$$
\lim_{\beta\to\infty}c(\beta)=\frac4\pi=1.2732\ldots,\qquad
\lim_{\beta\to0^+}c(\beta)=\frac{8}{\pi^2}=0.8105\ldots
$$
In particular, $c(\beta)$ would be discontinuous at $\beta=0$, since $c(0)=\sqrt6/\pi=0.7797\ldots<8/\pi^2$. This would be compatible with the change of limit shape at $\beta=0$, but the macroscopic limit shape does not by itself imply the required two-site limits.
\end{con}


\section*{Declaration of generative AI use}

During the preparation of this work, the authors used AI-assisted tools Codex and Claude to assist with language editing, manuscript organization, and exploration of some of the proof strategies and calculations. The authors reviewed and edited the resulting material and take full responsibility for the contents of the manuscript.

\providecommand{\bysame}{\leavevmode\hbox to3em{\hrulefill}\thinspace}
\providecommand{\MR}{\relax\ifhmode\unskip\space\fi MR }
\providecommand{\MRhref}[2]{%
  \href{http://www.ams.org/mathscinet-getitem?mr=#1}{#2}
}
\providecommand{\href}[2]{#2}

\end{document}